\documentclass[11pt,oneside,a4paper,reqno]{amsart}
\RequirePackage{fix-cm} 
\usepackage[T1]{fontenc}
\usepackage[english]{babel}
\usepackage{amssymb,amsmath,amsthm,graphicx,amsfonts}
\usepackage[foot]{amsaddr}
\usepackage{lmodern}

\usepackage[dvipsnames]{xcolor}
\definecolor{ourcolor}{RGB}{169,42,1}

\usepackage[
   rm={oldstyle=false,proportional=true},
   sf={oldstyle=true,proportional=true},
   tt={oldstyle=true,proportional=true,variable=true},
   qt=false,
 ]{cfr-lm}
\usepackage{cite,enumerate,setspace}
\usepackage[normalem]{ulem}
\usepackage{hyperref}
 \hypersetup{
 colorlinks,
 citecolor=ourcolor,
 linkcolor=ourcolor,
 urlcolor=black}

\numberwithin{equation}{section}

\usepackage{geometry}
\usepackage{caption}
\let\originalleft\left
\let\originalright\right
\renewcommand{\left}{\mathopen{}\mathclose\bgroup\originalleft}
\renewcommand{\right}{\aftergroup\egroup\originalright}

 \makeatletter
 \def\@seccntformat#1{\hspace*{0mm}%
  \protect\textup{\protect\@secnumfont
    \ifnum\pdfstrcmp{subsection}{#1}=0 \bfseries\fi
    \csname the#1\endcsname
    \protect\@secnumpunct
      }%
 }

\newcommand{\assign}{:=}
\newcommand{\mathLaplace}{\Delta}
\newcommand{\mathd}{\mathrm{d}}
\newcommand{\of}{:}
\newcommand{\tmabbr}[1]{#1}
\newcommand{\tmcolor}[2]{{\color{#1}{#2}}}
\newcommand{\tmem}[1]{{\em #1\/}}
\newcommand{\tmmathbf}[1]{\ensuremath{\boldsymbol{#1}}}
\newcommand{\tmname}[1]{\textsc{#1}}
\newcommand{\tmop}[1]{\ensuremath{\operatorname{#1}}}

\newcommand{\tmsamp}[1]{\textsf{#1}}

\newcommand{\tmtextup}[1]{\text{{\upshape{#1}}}}
\newenvironment{enumerateroman}{\begin{enumerate}[i.] }{\end{enumerate}}
\newtheorem{definition}{Definition}
{\theoremstyle{remark}\newtheorem{example}{Example}}
{\theoremstyle{remark}\newtheorem{remark}{Remark}}
\newtheorem{theorem}{Theorem}

\newcommand{\Om}{\Omega}

\newcommand{\JW}{W}
\newcommand{\jom}{\omega}
\newcommand{\ju}{\tmmathbf{u}}
\newcommand{\jm}{\tmmathbf{m}}
\newcommand{\jx}{x}
\newcommand{\jgam}{\gamma}
\newcommand{\jlam}{\lambda}
\newcommand{\Tsurf}{T}
\newcommand{\Ssurf}{S}
\newcommand{\ee}{\tmmathbf{e}}

\newcommand{\tfi}{\tmmathbf{\tau}_{\phi}}
\newcommand{\tti}{\tmmathbf{\tau}_t}
\newcommand{\Rm}{A}
\newcommand{\xip}{\tmmathbf{\xi}}

\newcommand{\cpt}{\tmtextup{)}}
\newcommand{\opt}{\tmtextup{(}}

\newcommand{\T}{\mathsf{T}}

\newcommand{\grad}{\nabla}
\newcommand{\divv}{\mathrm{div\,}}

\newcommand{\RR}{\mathbb{R}}

\newcommand{\Stwo}{\mathbb{S}}

\newcommand{\eqs}{=}

\let\oldtocsection=\tocsection

\let\oldtocsubsection=\tocsubsection

\let\oldtocsubsubsection=\tocsubsubsection

\renewcommand{\tocsection}[2]{\hspace{0em}\oldtocsection{#1}{#2}}
\renewcommand{\tocsubsection}[2]{\hspace{1em}\oldtocsubsection{#1}{#2}}
\renewcommand{\tocsubsubsection}[2]{\hspace{2em}\oldtocsubsubsection{#1}{#2}}

\begin{document}

\title[{\sc Minimizing harmonic maps
with axial symmetry}]{{\large \sc \bf \uppercase{Sufficient conditions for the existence of \\ minimizing harmonic maps
with axial symmetry \\ in the small-average regime}}}

\author[1]{{\sc Giovanni Di Fratta}$^{(1)}$}
\author[2]{{\sc Valeriy V. Slastikov}$^{(2)}$}
\author[3]{{\sc Arghir D. Zarnescu}$^{(3,4,5)}$}

\address[1]{ Dipartimento di Matematica e Applicazioni ``R.~Caccioppoli'', Università degli Studi di Napoli
    ``Federico II'', Via Cintia, 80126, Napoli, Italy.}
    
\address[2]{ School of Mathematics, University of
Bristol, Bristol BS8 1TW, United Kingdom.
}
    
\address[3]{ BCAM, Basque Center for Applied Mathematics, Mazarredo 14, E48009 Bilbao, Bizkaia, Spain.}   
\address[4]{ IKERBASQUE, Basque Foundation for Science, Maria Diaz de Haro 3, 48013, Bilbao, Bizkaia, Spain.} 
\address[5]{ Simion Stoilow Institute of Mathematics of the Romanian Academy, P.O. Box 1-764, RO-014700 Bucharest, Romania.}

\begin{abstract}
  The paper concerns the analysis of global minimizers of a Dirichlet-type
  energy functional defined on the space of vector fields $H^1 (
  \Ssurf, \Tsurf )$, where $\Ssurf$ and $\Tsurf$ are surfaces of
  revolution. The energy functional we consider is closely related to a
  reduced model in the variational theory of micromagnetism for the analysis
  of observable magnetization states in curved thin films. We show that axially symmetric minimizers always exist, and if the target surface $\Tsurf$ is never flat, then any coexisting minimizer must have line symmetry. Thus, the minimization problem reduces to the computation of an optimal one-dimensional profile. We also provide a necessary and sufficient condition for energy minimizers to be axially symmetric.
  
   \vspace{4pt}
  
  \noindent {\scriptsize {\sc Mathematics Subject Classification.} 35A23, 35R45, 49R05, 49S05, 82D40.}
  
  \vspace{1pt}
  
  \noindent {\scriptsize {\sc Key words.} Harmonic Maps, Axial Symmetry, Micromagnetics, Magnetic Skyrmions.}
\end{abstract}

\begingroup
\def\uppercasenonmath#1{} 
\let\MakeUppercase\relax 
\maketitle
\endgroup

\section{Introduction and motivation}
\noindent For given surfaces of revolution $\Ssurf, \Tsurf \subseteq \RR^3$ around the $\ee_3$ axis, we consider
the Dirichlet-type energy defined for every $\jm \in H^1( \Ssurf, \Tsurf)$ by
\begin{equation}
  \mathcal{E}_{\jom} ( \jm ) \assign \int_{\Ssurf} | \grad \jm
  (\xi) |^2 \mathd \xi + \int_{\Ssurf} g ( \jm (\xi) \cdot
  \tmmathbf{a} (\xi) ) \mathd \xi + \int_{\Ssurf} | \langle
  \jm \rangle_{\Ssurf_{\xi}} \times \ee_3 |^2  \jom^2 (\xi) \mathd
  \xi, \label{eq:mainenfunc}
\end{equation}
where $g : \RR \to \RR_+$ is an anisotropic potential, $\tmmathbf{a}:
\Ssurf \to \RR^3$ is a prescribed vector field that can be either
axially symmetric or axially antisymmetric (see
section~\ref{subsubsec:symvecsfields}), $\langle \jm
\rangle_{\Ssurf_{\xi}}$ is the average of $\jm$ along the circle of
latitude $\Ssurf_{\xi} \assign ( \xi \cdot \ee_3 ) \ee_3 + \Stwo^1
( \xi \cdot \ee_1 )$, i.e.,
\begin{equation}
  \langle \jm \rangle_{\Ssurf_{\xi}} = \frac{1}{| \Ssurf_{\xi}
  |} {\int_{\Ssurf_{\xi}}}  \jm (\xi) \mathd \xi,
\end{equation}
and $\jom : \Ssurf \to \RR_+$ is a generic measurable function that
weights the strength of the tendency of $\langle \jm
\rangle_{\Ssurf_{\xi}}$ to be aligned along the unit vector $\ee_3$.
Detailed hypotheses on the regularity of the involved quantities, as well as
precise definitions of the terms employed, will be given in
section~\ref{sec:setup}.

The main aim of this paper is to show that under mild conditions on the weight
function $\jom$, every global minimizer of $\mathcal{E}_{\jom}$ has axial
symmetry. When $\tmmathbf{a}$ is axially symmetric, a particular case of our findings gives a necessary and sufficient condition for energy minimizers to be axially symmetric: the existence of axially symmetric energy minimizers is equivalent to the existence of \emph{axially null-average} minimizers, i.e., to the existence of minimizers $\jm \in H^1( \Ssurf, \Tsurf)$ such that $\langle \jm \rangle_{\Ssurf_{\xi}} \times \ee_3 = 0$ for every $\xi \in \Ssurf$ (see section~\ref{sec:furthresappls} and Theorem~\ref{thm:main3} in there). This characterization also holds when the last term in \eqref{eq:mainenfunc} is absent. Note that while it is always the case that axially symmetric configurations are axially null average, minimality allows for the converse implication.

In what follows, to shorten notation and enhance readability, we set
\begin{equation}
  \mathcal{D}_{\Ssurf} ( \jm ) \assign \int_{\Ssurf} | \grad
  \jm (\xi) |^2 \mathd \xi, \quad \mathcal{A}_{\Ssurf} ( \jm
  ) \assign \int_{\Ssurf} g ( \jm (\xi) \cdot \tmmathbf{a} (\xi)
  ) \mathd \xi, \label{eq:energyterms0}
\end{equation}
and
\begin{equation}
  \mathcal{P}_{\Ssurf} ( \jm ) \assign \int_{\Ssurf}|
  \langle \jm \rangle_{\Ssurf_{\xi}} \times \ee_3 |^2  \jom^2
  (\xi) \mathd \xi, \label{eq:energyterms1}
\end{equation}
so that the energy functional we are interested in takes the form
\begin{equation}
  \mathcal{E}_{\jom} ( \jm ) =\mathcal{D}_{\Ssurf} ( \jm
  ) +\mathcal{A}_{\Ssurf} ( \jm ) +\mathcal{P}_{\Ssurf}
  ( \jm ) .
\end{equation}
Our investigations fit with that thread of results concerned with the
study of harmonic maps with symmetries. Indeed, when $\Ssurf$ is a surface
with boundary, our analysis applies to the existence of axially symmetric
solutions of  the Euler-Lagrange equations of $\mathcal{E}_{\jom}$, which are harmonic-type equations between surfaces of revolution.

Also, when $\Tsurf \assign \Stwo^2$, the energy functional
$\mathcal{E}_{\jom}$ is closely related to a reduced model in the variational
theory of micromagnetism, a model  relevant for the description of observable magnetization states
in curved thin films. As explained in section~\ref{subsec:physcontext}, the functional
$\mathcal{E}_{\jom}$ can be considered a model of the micromagnetic energy functional in the asymptotic regime of curved thin films where, in addition to magnetocrystalline and shape anisotropies,
higher-order magnetostatic effects are taken into account through the term $\mathcal{P}_{\Ssurf}$, which favors configurations that are
$\Ssurf_{\xi}$-null-average in the perpendicular component. According to our
findings, micromagnetic ground states always have axial symmetry in the curved
thin-film regime and under the influence of nonnegligible higher-order
magnetostatic effects. The conclusions we reach can be applied to other
significant physical systems, e.g., to understand the existence of symmetric
textures in the Oseen-Frank theory of nematic liquid
crystals~{\cite{Golovaty_2017,Nitschke_2018}}.

\subsection{Outline}In the next section~\ref{subsec:physcontext}, we briefly
present the physical context that led us to the investigation of
\eqref{eq:mainenfunc}, so to give the reader a broader perspective on the
relevance of the energy functional $\mathcal{E}_{\jom}$. In
section~\ref{subsec:stateofart}, after a brief review of earlier research on
the symmetry of harmonic maps, we discuss previous studies on the symmetry
properties of minimizers of the micromagnetic energy functional in the curved
thin film regime. In section~\ref{sec:setup}, we describe the rigorous setting
of the problem and detail the contributions of the present work. Proofs are
given in section~\ref{sec:proofs}. In section~\ref{sec:furthresappls}, we
present further results and some applications to the existence of axially
symmetric solutions of elliptic PDEs.

\subsection{Physical context}\label{subsec:physcontext}Over the last decade,
considerable experimental and theoretical research has been done on the
physics of ferromagnetic systems with curved shapes. One of the main reasons
is that the curved geometry can lead to effective antisymmetric interactions
and, as a result, to the consequent formation of magnetic skyrmions, i.e.,
chiral spin textures with a non-trivial topological degree, even in the
absence of spin-orbit coupling mechanism, responsible for
Dzyaloshinskii--Moriya interactions (DMI). The evidence of these states sheds
light on the role of geometry in magnetism: magnetic skyrmions can be
stabilized by curvature effects only, in contrast to the planar case where DMI
is required~{\cite{MaKSheka2022}}.

Also, recent advances in the nanofabrication of magnetic hollow particles
have sparked interest in these geometries, which lead to artificial materials
with unexpected characteristics and diverse applications spanning from logic
devices to biomedicine~{\cite{Streubel_2016}}.

From a mathematical point of view, the analysis of mesoscale magnetism in
curved geometry can be carried out within the framework of the variational
theory of micromagnetism {\cite{BrownB1963,Hubert1998,Di_Fratta_2019_var}},
where the order parameter is the magnetization $\jm$ subject to the
{\tmem{saturation}} constraint of being $\Stwo^2$-valued. In this framework,
the energy term $\mathcal{A}_{\Ssurf}$ accounts for the so-called
{\tmem{crystal}} and {\tmem{shape-anisotropy effects}}. Indeed, it is
well-known that when the ferromagnet occupies a thin shell whose thickness is
significantly smaller than the size of the system, the dominant energy
contribution is encapsulated in the energy functional
{\cite{GioiaJames97,carbou2001thin,Di_Fratta_2020}}
\begin{equation}
  \mathcal{F}_{\kappa} : \jm \in H^1( \Ssurf, \Stwo^2) \mapsto
  \int_{\Ssurf} | \grad \jm (\xi) |^2 \mathd \xi + \kappa
  \int_{\Ssurf} ( \jm (\xi) \cdot \tmmathbf{n} (\xi) )^2 \mathd
  \xi, \label{eq:sphericshelllimitgen}
\end{equation}
where $\Ssurf$ is the surface generating the hollow nanoparticle by extrusion
along the normal direction $\tmmathbf{n}$, and $\kappa \in \RR$ is an
effective anisotropy parameter accounting for both {\tmem{shape}} and
{\tmem{crystal}} anisotropy. For large $\kappa > 0$, tangential vector fields
are energetically favored; for large $\kappa < 0$, i.e., when perpendicular crystal anisotropy
prevails over shape anisotropy, energy minimization prefers
normal vector fields. Note that, from the variational point of view, the
expression of $\mathcal{F}_{\kappa}$ is nothing but a particular case of the
energy functional $\mathcal{E}_{\jom}$ when $\Tsurf \assign \Stwo^2$,
$\tmmathbf{a} (\xi) \assign \tmmathbf{n}_{\Ssurf} (\xi) $ and $g (s) \assign
\kappa s^2$ if $\kappa > 0$ and $g (s) \assign | \kappa | (1 - s^2)$ if
$\kappa < 0$. In other words, our energy term $\mathcal{A}_{\Ssurf}$ in
\eqref{eq:energyterms0} accounts for possible generalizations of the second
term in \eqref{eq:sphericshelllimitgen}. Other typical expressions of
$\tmmathbf{a}$ and $g$ issuing from the variational theories of micromagnetics
and nematic liquid crystals are $\tmmathbf{a} (\xi) \assign
\tmmathbf{n}_{\Ssurf} (\xi) $, $\tmmathbf{a} (\xi) \assign
\tmmathbf{n}_{\Tsurf} (\xi)$, $\tmmathbf{a} (\xi) = \ee_3$, or $g (s) =
\lambda (1 - s^2)^2$ for some $\lambda \in \RR_+$.

In
the language of the variational theory of micromagnetism, we can think of the energy term $\mathcal{P}_{\Ssurf}$ as accounting for
higher-order effects in the long-range magnetic dipole-dipole interactions. Indeed, it is well-known that
in the classical three-dimensional setting, the magnetostatic self-energy
associated with a distribution of magnetization $\jm \in H^1( \Om, \RR^3)$ in a domain $\Om \subseteq \RR^3$, favors solenoidal vector fields,
i.e., divergence-free vector fields that are tangential to the
boundary~{\cite{BrownB1962,BrownB1963}}. From the variational point of view,
this means that the magnetostatic self-energy is minimized when
\begin{equation}
  \divv  \jm = 0 \quad \text{in } \Om, \quad \jm \cdot \tmmathbf{n}_{\partial
  \Om} = 0 \quad \text{on } \partial \Om . \label{eq:nullaveragefavor}
\end{equation}
Now, as shown in {\cite{GioiaJames97,carbou2001thin,Di_Fratta_2020}}, but also
apparent from \eqref{eq:sphericshelllimitgen}, at the leading order in the
energy reduction from $3 d$ to $2 d$, only the tendency of being tangential to
the boundary is preserved. At the same time, any aspect of the divergence-free
conditions is lost. To explain this last point better, we recall that any
vector-field $\jm \in H^1 ( \Om, \RR^3 )$ satisfying
\eqref{eq:nullaveragefavor} is necessarily null-average in $\Om$. Indeed, for
$i \in \{ 1, 2, 3 \}$, we have $\divv ( \jm x_i ) = x_i \divv \jm +
\jm \cdot \ee_i$; therefore, if $\jm$ satisfies \eqref{eq:nullaveragefavor},
then $\langle \jm \cdot \ee_i \rangle_{\Om} = \langle \divv
( \jm x_i ) \rangle_{\Om} = 0$, because of the divergence
theorem and the tangential boundary condition on $\jm$. Despite being a simple
mathematical consequence of \eqref{eq:nullaveragefavor}, the favoring of
null-average configurations is often the main qualitative property used to
describe the physical effect of the demagnetizing field (see,
e.g.,~{\cite{BertottiB1998}}).

Maintaining some of these structural consequences of
\eqref{eq:nullaveragefavor} in the passage from $3 d$ to $2 d$ would be
desirable, and this is the main motivation for considering the
term $\mathcal{P}_{\Ssurf}$ in our analysis: to keep track of some of the
qualitative features of the demagnetizing field that get lost at the leading
order in the energy expansion. Indeed, if $\Om \assign \Sigma \times [0, 1]$,
$\Sigma \subseteq \RR^2$, is a cylindrical domain, and $\jm$ satisfies
\eqref{eq:nullaveragefavor}, then other than being $\langle \jm
\rangle_{\Om} = 0$ we can say more, i.e., that for every $h \in [0, 1]$
there holds $\langle \jm \times \ee_3 \rangle_{\Sigma_h} = 0$, with
$\Sigma_h \assign \Sigma \times \{ h \}$. In fact,
\[ \int_{{\Sigma_h} } \jm (x) \cdot \ee_i \mathd x \eqs \int_{{\partial
   \Sigma_h} } x_i  \jm (\sigma) \cdot \tmmathbf{\nu}_{{\partial \Sigma_h} }
   (\sigma) \mathd \sigma = 0 \]
because $\tmmathbf{\nu}_{{\partial \Sigma_h} } (\sigma)
=\tmmathbf{n}_{\partial \Om} (\sigma)$ for every ${\sigma \in \partial
\Sigma_h} $. Although the argument fails for noncylindrical domains, and it is
not clear to what extent the tendency towards $\langle \jm \times \ee_3
\rangle_{\Sigma_h} = 0$ persists in general surfaces of revolution, what
we pictured motivates us enough to explore the impact that a term
like $\mathcal{P}_{\Ssurf}$ in \eqref{eq:energyterms1} has on the energetic
landscape.

\subsection{State of the art}\label{subsec:stateofart}The minimization problem
for Dirichlet-type energy functionals between manifolds naturally arises in
differential geometry in the study of unit-speed geodesics and minimal
submanifolds~{\cite{Eells1995,Schoen1997}}. It also occurs in nonlinear field
theories every time local interactions are brought into the realm of elastic
deformations. Prominent examples are the elastic free energy in the
Oseen-Frank theory of nematic liquid crystals {\cite{Ball_2017,Virga_2018}},
symmetric and antisymmetric exchange interactions in the variational theory of
micromagnetism~{\cite{BrownB1963,Hubert1998}}, and
$M$-theory~{\cite{Albeverio_1997,Becker_2006}}. The systematic treatment
gained impetus after the seminal work of Eells and Sampson {\cite{Eells1964}},
who showed that under certain technical conditions on the base and target
manifolds, every continuous map is homotopic to a harmonic map.

The existence of minimizing harmonic maps with symmetry has been the topic of
several studies, typically set around specific geometries. In
{\cite{Sandier_1993,sandier1993symmetry}}, it is shown that for any symmetrical domain in $\RR^2$
and any symmetrical boundary datum that takes values in a closed hemisphere,
minimizing harmonic maps must be radially symmetric. In {\cite{MR0984634}},
using projection and averaging procedures introduced in~{\cite{Coron1989}},
Coron and Helein study harmonic diffeomorphisms between the Euclidean $n$-ball
(minus a finite number of points) and a Riemannian manifold. They show that
under some conditions on the involved parameters, there always exists a smooth
$\tmop{SO} (3)$-equivariant map defined on $B^3$ with values in $E^3 (a)
\assign \{(x, y) \in \RR^3 \times \RR \of \hspace{0.22em} |x|^2 + y^2 / a^2 =
1\}$ improving some earlier results of Baldes~{\cite{Baldes_1984}}. In
{\cite{Hardt_1992}}, \ Hardt, Lin, and Poon investigated the existence of
$\Stwo^2$-valued axially symmetric harmonic maps with prescribed singularities
(see also {\cite{MR1178093,Hardt_1992}} and {\cite{Martinazzi_2011}}).

The literature on the subject is vast, and only a systematic review could give
due recognition to the results obtained over the years. Some references are available in \cite{brezis1999symmetry, pisante2014symmetry}. Here, we limit
ourselves to the literature intimately related to our investigations on the
energy functional $\mathcal{E}_{\jom}$, and we refer the reader to
{\cite{Eells1995,Schoen1997,H_lein_2002,Lin2008}} \ and the references therein
for further results on the topic.

Magnetic hollow nanoparticles usually are of the shape of a surface of
revolution, and symmetry properties of the ground-states are usually derived
in specific geometries. For spherical thin films, i.e., when $\Ssurf =
\Stwo^2$, we showed in {\cite{Di_Fratta_2019}} that when $\kappa \leqslant -
4$, the normal vector fields $\pm \tmmathbf{n}$ are the {\tmem{only}}
{\tmem{global}} minimizers of the energy functional  $\mathcal{F}_{\kappa}$ defined in \eqref{eq:sphericshelllimitgen}. The interest in results
of this kind is in the topological remark that $\pm \tmmathbf{n}$ carry
different {\tmem{skyrmion numbers}} because of deg$(\pm \tmmathbf{n}) = \pm 1$
---see also {\cite{Melcher_2019}} for stationary states topologically distinct
from the ground state. In physical terms, this translates into their
robustness against thermal fluctuations and external perturbations with
far-reaching consequences for modern magnetic storage
technologies~{\cite{Fert_2013}}. The vector fields $\pm \tmmathbf{n}$ have
full rotational symmetry, and their local stability persists up to $\kappa <
0$. They lose stability for $\kappa > 0$ when, therefore, new ground states
have to appear. No more than this is currently known. The reason is that when
$\kappa > - 4$, the energy landscape of $\mathcal{F}_{\kappa}$ is very hard to
describe analytically, and it is unclear which aspects of the symmetry of the
problem are retained in the shape of minimizers. Numerics suggest that when
$\kappa > 0$, the energy $\mathcal{F}_{\kappa}$ exhibits magnetic states with
skyrmion numbers $0$, $\pm 1$, all having axial
symmetry~{\cite{Kravchuk_2012,Kravchuk_2016}}. However, no mathematical
evidence of this has been established. Our paper aims to gain some insight
into the question.

In {\cite{Fratta_2023}}, it
is shown that when $\Ssurf = \Stwo^1 \times [0, 1]$, the $\Stwo^1$-valued
normal vector fields $\pm \tmmathbf{n}$ are the unique ground state when
provided that $\kappa$ is sufficiently negative ($\kappa \leqslant - 3$), but
almost nothing is currently known when $\kappa > - 3$. The only thing one can
say is that if the ground states are, as we expect, in-plane, then global
minimizers lose axial symmetry because they are no more null average
(see~{\cite[Figure~5]{Fratta_2023}} and {\cite{Carbou_2022}}).

However, the axial symmetry of the minimizers certainly fails, for example,
when $\Ssurf = \mathbb{D}$ is the unit disk in $\RR^2$. Indeed, on the one
hand, as we show in {\cite{Di_Fratta_2022}}, the absence of curvature favors
radial ground states rather than axially symmetric, i.e., minimizers satisfy
the relation $\jm (x) = \jm (| x |)$ for every $x \in \mathbb{D}$. On the
other hand, our proofs still work if one adds to the expression of
$\mathcal{E}_{\jom}$ a penalization on the boundary as in
{\cite{Di_Fratta_2022}}. Therefore, we conclude that in the generality we are
treating the problem, it is essential to consider the term
$\mathcal{P}_{\Ssurf}$ to obtain the existence of axially symmetric minimizers
regardless of the specific choices of $\Ssurf$ and $\Tsurf$ in the class of
surfaces of revolution.

\section{Contributions of the present work}\label{sec:setup}

\subsection{Notation and setup}The main results of the present work concern
the symmetry properties of energy minimizing maps. To state our results
precisely, we need to set up the framework, the mathematical notation, and the
terminology used throughout the paper.

\subsubsection{Sobolev spaces on surfaces}For a given $C^1$-surface $\Ssurf
\subset \RR^3$, we denote by $H^1 ( \Ssurf, \RR^3 )$ the Sobolev
space of vector-valued functions defined on $\Ssurf$ endowed with the norm
\begin{equation}
  \| \jm \|^2_{H^1( \Ssurf, \RR^3 )} \assign
  \int_{\Ssurf} | \jm (\xi) |^2 \mathd \xi + \int_{\Ssurf} |
  \grad \jm (\xi) |^2 \mathd \xi .
\end{equation}
Here, $\grad$ is the tangential gradient of $\jm$ at $\xi \in \Ssurf$, and
$| \grad \jm (\xi) |^2 = \sum_{i = 1}^2 |
\partial_{\tmmathbf{\mu}_i (\xi)} \jm (\xi) |^2$ if $(\tmmathbf{\mu}_1
(\xi), \tmmathbf{\mu}_2 (\xi))$ is an orthonormal basis of $T_{\xi} \Ssurf$.
Also, given two surfaces $\Ssurf, \Tsurf \subseteq \RR^3$, we write $H^1
( \Ssurf, \Tsurf )$ for the metric subspace of $H^1( \Ssurf,
\RR^3 )$ made by vector-valued functions with values in $\Tsurf$.

\subsubsection{Surfaces of revolution}In what follows we denote by $I
\subseteq \RR$ a closed interval, by $( \ee_1, \ee_2, \ee_3 )$ the
standard ordered basis of $\RR^3$, and by $\Rm^{\top} (\phi)$ the rotation
matrix about the $\ee_3$-axis given by
\begin{equation}
  \Rm^{\T} (\phi) \assign \left(\begin{array}{ccc}
    \cos \phi & - \sin \phi & 0\\
    \sin \phi & \cos \phi & 0\\
    0 & 0 & 1
  \end{array} \right) . \label{ATMatrix}
\end{equation}
By a regular simple curve, we mean the image of a $C^1$-map $\jgam \of I
\mapsto \RR^3$ such that $\dot{\jgam} (t) \neq 0$ for every $t \in I$, and
with no self-intersections, i.e., such that the only possible loss of
injectivity in $\jgam$ arises at the endpoints of $I$, case in which the curve
closes into a loop. As customary, we often refer to $\jgam$ as a curve rather
than just its image.

\begin{figure}[t]
  {\includegraphics[width=\linewidth]{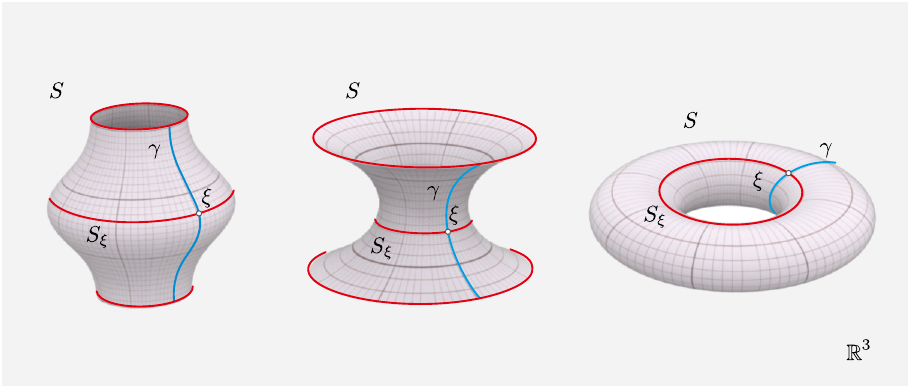}}
  \caption{\label{fig:surfsrevol}Given a regular curve $\jgam$ in the $\jx,
  z$-plane, which lies at a nonnegative distance from the $\ee_3$-axis, the
  surface of revolution $\Ssurf$ generated by $\jgam$ is the set $\Ssurf
  \assign \cup_{\xi \in \jgam} \Ssurf_{\xi}$, where $\Ssurf_{\xi} \assign
  ( \xi \cdot \ee_3 ) \ee_3 + \Stwo^1 ( \xi \cdot \ee_1
  )$ is the circle of latitude at $\xi \in \Ssurf$. The resulting
  surface can have a boundary or not.}
\end{figure}

Given a regular simple curve $\jgam \of t \in I \mapsto \jgam (t) = ( \jx
(t), 0, z (t) ) \in \RR^3$, with $\jx (t) \geqslant 0$, the surface of
revolution $\Ssurf \subseteq \RR^3$ generated by $\jgam$ is the image of the
parameterization defined for $0 \leqslant \phi \leqslant 2 \pi$ and $t \in I$
by
\begin{equation}
  \xip (\phi, t) = \Rm^{\top} (\phi) \jgam (t) . \label{eq:paramS}
\end{equation}
In more intrinsic terms, given a regular simple curve $\jgam$ in the $\jx,
z$-plane, which lies at a nonnegative distance from the $\ee_3$-axis, the
surface of revolution $\Ssurf$ generated by $\jgam$ is the set $\Ssurf \assign
\cup_{\xi \in \jgam} ( ( \xi \cdot \ee_3 ) \ee_3 + \Stwo^1
( \xi \cdot \ee_1 ) )$, where $\Stwo^1 (r)$ is the circle in
$\RR^2 \times \{ 0 \}$ centered at the origin and of radius $r > 0$. For our
purposes, it is convenient to denote by $\Ssurf_{\xi} \assign ( \xi \cdot
\ee_3 ) \ee_3 + \Stwo^1 ( \xi \cdot \ee_1 )$ the circle of
latitude at $\xi \in \Ssurf$. After that, we have that
\begin{equation}
  \Ssurf \assign \cup_{\xi \in \jgam} \Ssurf_{\xi} .
  \label{eq:Sasunioncircleslat}
\end{equation}
Given that $\jgam$ is defined on a closed interval, the resulting surface of
revolution is always topologically closed and possibly with a boundary (see
Figure~\ref{fig:surfsrevol}).

The tangent space to $\Ssurf$ at $\xip (\phi, t)$ is generated by the two
vectors
\begin{align}
  \tfi (\phi, t) \assign \partial_{\phi} \xip (\phi, t) & \eqs 
  \partial_{\phi} \Rm^{\top} (\phi) \jgam (t),  \label{eq:tauphi}\\
  \tti (\phi, t) \assign \partial_t \xip (\phi, t) & \eqs  \Rm^{\top} (\phi)
  \dot{\jgam} (t) .  \label{eq:taut}
\end{align}
Given that $\Rm (\phi) \partial_{\phi} \Rm^{\top} (\phi) v = \ee_3 \times v$
for every $v \in \RR^3$, one easily gets that
\begin{equation}
  \tfi (\phi, t) \cdot \tti (\phi, t) = \Rm (\phi) \partial_{\phi} \Rm^{\top}
  (\phi) \jgam (t) \cdot \dot{\jgam} (t) = ( \ee_3 \times \jgam (t)
  ) \cdot \dot{\jgam} (t) = ( \jgam (t) \times \dot{\jgam} (t)
  ) \cdot \ee_3 = 0,
\end{equation}
because $\jgam (t) \times \dot{\jgam} (t)$ is directed along $\ee_2$. Thus,
the coordinate system is {\tmem{orthogonal}} but, in general, {\tmem{not}}
orthonormal. Indeed, the metric coefficients are given by
\begin{equation}
  \mathfrak{h}_1 (t) \assign | \tfi | = | \jx (t) |,
  \quad \mathfrak{h}_2 (t) \assign | \tti | = | \dot{\jgam}
  (t) | . \label{eq:metrich1}
\end{equation}
Since $\tfi (\phi, t)$ and $\tti (\phi, t)$ are orthogonal, the area element
on $\Ssurf$ assumes the form
\begin{equation}
  \sqrt{\mathfrak{g} (t)} \eqs | \partial_{\phi} \xip \times \partial_t
  \xip | = | \tfi (\phi, t) | \cdot | \tti (\phi, t)
  | =\mathfrak{h}_1 (t) \cdot \mathfrak{h}_2 (t) . \label{eq:areaeleS}
\end{equation}
Note that the area element in \eqref{eq:areaeleS} depends only on the $t$
variable and, by the regularity hypotheses in Remark~\ref{rmk:regulsurf}, we
have that $\sqrt{\mathfrak{g} (t)} > 0$ for every $t \in I$ with the possible
exception of the two boundary points of $\jgam$ if $\jgam$ touches the
$\ee_3$-axis.

\begin{remark}
  \label{rmk:regulsurf}({\tmem{On the regularity of surfaces of revolution}})
  Throughout the paper, we always assume that the generated surface of
  revolution $\Ssurf$ is {\tmem{regular}}, i.e., that the parameterization
  \eqref{eq:paramS} is of class $C^1$, and that the induced area element
  \eqref{eq:areaeleS} never vanishes. Now, since the generating curve $\jgam$
  is regular, the surface of revolution $\Ssurf$ induced by $\jgam$ is
  certainly regular whenever $\varrho_{\jgam} \assign \inf_{t \in I}  \jx (t)
  > 0$. Instead, if $\varrho_{\jgam} = 0$, we assume, as in the case of the
  meridian arc generating $\Stwo^2$, that $\jgam$ touches the $\ee_3$-axis
  perpendicularly and (at most) at two distinct points so that a smooth
  surface of revolution results. In what follows, whenever we name a surface
  of revolution, it is always understood that the previous regularity
  assumptions are met.
\end{remark}

\subsubsection{Symmetric vector fields}\label{subsubsec:symvecsfields}The
results of our paper guarantee the existence of minimizers with specific
symmetry properties. The notion of an axially symmetric vector field is standard.

\begin{definition}
  \label{def:axialsymmetry}With $\Ssurf, \Tsurf \subseteq \RR^3$ being two
  surfaces of revolution, we say that a vector field $\ju : \Ssurf \to
  \Tsurf$ is {\tmem{axially symmetric}} when the following property holds:
  \begin{equation}
    \ju ( \Rm^{\top} (\phi) \xi ) \eqs \Rm^{\top} (\phi)
    \ju (\xi) \quad \forall \phi \in \RR, \forall \xi \in \Ssurf .
    \label{eq:coordinatefreeaxsymm}
  \end{equation}
  We say that $\ju$ is {\tmem{axially antisymmetric}} if
  \begin{equation}
    \ju ( \Rm^{\top} (\phi) \xi ) \eqs \Rm (\phi) \ju
    (\xi) \quad \forall \phi \in \RR, \forall \xi \in \Ssurf .
    \label{eq:coordinatefreeaxsymmanti}
  \end{equation}
  We say that a vector field $\ju$ has {\tmem{line symmetry}} if for every
  fixed $\xi\in \Ssurf$ one has either $\ju ( \Rm^{\top} (\phi)
  \xi) = \Rm^{\top} (\phi) \ju (\xi)$ for every $\phi \in \RR$ or $\ju( \Rm^{\top} (\phi)
  \xi ) = \Rm (\phi) \ju (\xi)$ for every $\phi \in \RR$.
\end{definition}

\begin{remark}
  Surfaces of revolution are $G$-sets with respect to the action of the group
  $G$ of rotations around the $\ee_3$-axis. In the language of transformation
  groups, an axially symmetric vector field can be equivalently defined as an
  equivariant map with respect to the $G$-sets $\Ssurf$ and $\Tsurf$, where
  $G$ is the group of rotations around the $\ee_3$-axis. Similarly, one could
  define an axially {\tmem{anti}}symmetric vector field as an
  equi{\tmem{contra}}variant map with respect to the $G$-sets $\Ssurf$ and
  $\Tsurf$.
\end{remark}

In local coordinates, and with the convenient abuse
of notation $\ju (\phi, t) \assign ( \ju \circ \xip ) (\phi, t)$,
the vector field $\ju : \Ssurf \to \Tsurf$ is axially symmetric if,
and only if, $\ju (\phi, t) = \Rm^{\top} (\phi) \ju (0, t)$ for every $(\phi,
t) \in \RR \times I$, and is axially antisymmetric if, and only if, $\ju
(\phi, t) = \Rm (\phi) \ju (0, t)$ for every $(\phi, t) \in \RR \times I$.
Instead, in general, if $\ju$ has line symmetry, then we are still allowing
for the existence of $t_1 \neq t_2 \in I$ such that $\ju (\phi, t_1) =
\Rm^{\T} (\phi) \ju (0, t_1)$ and $\ju (\phi, t_2) = \Rm (\phi) \ju (0, t_2)$
for every $\phi \in \RR$.


Note that the class of axially symmetric vector-field is essentially disjoint
from the class of axially antisymmetric vector fields. Indeed, if $A^{\T}
(\phi) \tmmathbf{\alpha} (t) = A (\phi) \tmmathbf{\beta} (t)$ for
$\Tsurf$-valued profiles $\tmmathbf{\alpha}$ and $\tmmathbf{\beta}$, then
necessarily $\tmmathbf{\alpha}=\tmmathbf{\beta}$ and $\tmmathbf{\alpha}$ is
directed along the $\ee_3$-axis. Also, note that the vector field
$\tmmathbf{v} (\phi, t) = A (\phi) \tmmathbf{\alpha} (t)$ is axially
antisymmetric if, and only if, $\ju (\phi, t) =\tmmathbf{v} (- \phi, t)$ is
axially symmetric.

We provide some examples to help clarify the notion of line symmetry. Suppose $\gamma : I \to \RR^3$ is the
image of the curve generating the surface of revolution $S$ and
$\tmmathbf{\alpha}: I \mapsto T$ is a smooth vector field along $\gamma$. The
vector field $\tmmathbf{u}: S \to T$ defined in local coordinates by
$\tmmathbf{u} (\phi, t) = A^{\top} (\phi) \tmmathbf{\alpha} (t)$ is axially
symmetric. The vector field $\tmmathbf{v}: S \to T$ defined in local
coordinates by $\tmmathbf{v} (\phi, t) = A (\phi) \tmmathbf{\alpha} (t)$ is
axially antisymmetric. Axially symmetric and antisymmetric vector fields
belong to the class of vector fields with line symmetries. However, there are examples of vector fields with line symmetry that are neither
axially symmetric nor axially antisymmetric. For that, consider the surface of
revolution $S = 2\mathbb{D}$ given by the disk of radius two centered at the
origin. If $\tmmathbf{v}: \mathbb{D} \to T$ is any axially
antisymmetric vector field such that $\tmmathbf{v}_{|
\partial \mathbb{D}} \equiv \tmmathbf{e}_3$, and $\tmmathbf{u}
\of 2\mathbb{D}\backslash\mathbb{D} \to T$ is any axially symmetric
vector field such that $\tmmathbf{u}_{| \partial \mathbb{D}}
\equiv \tmmathbf{e}_3$, then the vector field $\tmmathbf{w}: 2\mathbb{D}
\to T$ obtained by gluing $\tmmathbf{u}$ and $\tmmathbf{v}$ along $\partial \mathbb{D}$ is a
continuous vector field with line symmetry, which is neither axially symmetric
nor axially antisymmetric; the example is prototypical because it can be easily
generalized to build smooth line-symmetric vector fields defined on arbitrary surfaces of revolution that are not
axially symmetric nor axially antisymmetric.

\begin{figure}[t]
  {\includegraphics[width=\linewidth]{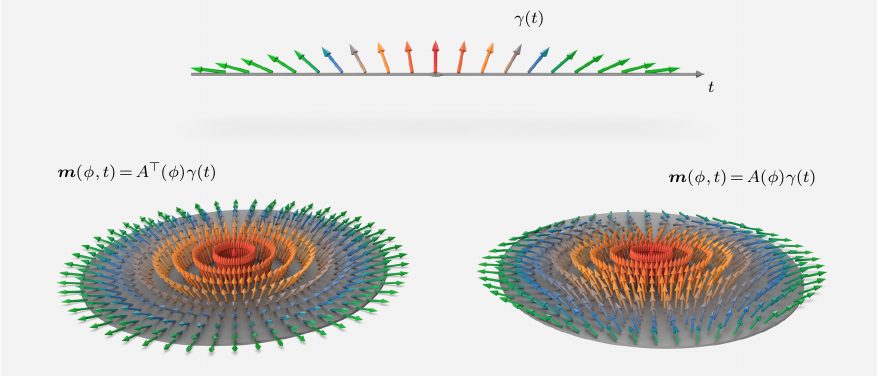}}
  \caption{\label{fig:diffasaas}Given the profile $\gamma : t \in I \mapsto
  (\sin  t, 0, \cos  t) \in \Stwo^2$, $I = [- 1, 1]$, the axially symmetric
  vector field $\jm (\phi, t) \assign A^{\T} (\phi) \gamma (t)$ (depicted on
  the left) and the axially antisymmetric vector field $\jm (\phi, t) \assign
  A (\phi) \gamma (t)$ (depicted on the right) can look pretty different.
  Nevertheless, they have the same Dirichlet energy.}
\end{figure}

We stress that for a given profile $\tmmathbf{\alpha}: I \to \Tsurf$,
the axially symmetric vector field $\Rm^{\top} (\phi) \tmmathbf{\alpha} (t)$
and the axially antisymmetric vector field $\Rm (\phi) \tmmathbf{\alpha} (t)$
generated by $\tmmathbf{\alpha}$ are, in general, different from the global
point of view. For example ({\tmabbr{cf.}}~Figure~\ref{fig:diffasaas}), for
$\tmmathbf{\alpha}: t \in I \mapsto (\sin  t, 0, \cos  t) \in \Stwo^2$, $I =
[- 1, 1]$, the axially symmetric vector field $A^{\T} (\phi) \tmmathbf{\alpha}
(t)$ has mirror symmetry with respect to every plane orthogonal to $\Stwo^1
\times \{ 0 \}$, while the generated axially antisymmetric vector field has
only a finite number of mirror symmetries.
However, note that both share the same Dirichlet energy
$\mathcal{D}_{\Ssurf}$ and the same energy
$\mathcal{P}_{\Ssurf}$. In our setting, the anisotropy term
$\mathcal{A}_{\Ssurf}$ is the only contribution that, from the variational
point of view, can ultimately favor one over the other.

\subsection{Statement of main results}\tmcolor{blue}{}Even though our primary
motivation comes from curved thin-film structures, we formulate and prove our
results for the case where the target space is a generic surface of revolution
rather than just $\Stwo^2$. With the usual abuse of notation, we set $\jom
(\phi, t) \assign \jom ( \xip (\phi, t) )$, $(\phi, t) \in [0, 2
\pi] \times I$, and, for every $t \in I$ we introduce the {\tmem{circular
integral weight}}
\begin{equation}
  \JW^2 (t) \assign \int_0^{2 \pi} \jom^2 (\phi, t) \mathd \phi .
  \label{eq:Womega}
\end{equation}
We recall that the metric coefficient $\mathfrak{h}_1 (t) \assign | \jx
(t) |$ measures the distance of the curve $\jgam$ generating $\Ssurf$
from the $\ee_3$-axis or, equivalently, the radius of the circle of latitude
$\Ssurf_{\jgam (t)}$. The first result of our paper is stated in the following
theorem.

\begin{theorem}
  \label{thm:main0} For given surfaces of revolution $\Ssurf, \Tsurf \subseteq
  \RR^3$, we consider the energy functional defined for every $\jm \in H^1
  ( \Ssurf, \Tsurf )$ by
  {\opt}cf.{\tmem{~\eqref{eq:mainenfunc}}}{\cpt}
  \[ \mathcal{E}_{\jom} ( \jm ) \assign \int_{\Ssurf} | \grad
     \jm (\xi) |^2 \mathd \xi + \int_{\Ssurf} g ( \jm (\xi) \cdot
     \tmmathbf{a} (\xi) ) \mathd \xi + \int_{\Ssurf}| \langle
     \jm \rangle_{\Ssurf_{\xi}} \times \ee_3 |^2  \jom^2 (\xi)
     \mathd \xi \]
  with $\jom : \Ssurf \to \RR_+$ a measurable weight such that $\sup_{t\in I} 
  (\mathfrak{h}_1 (t) \JW (t))<+\infty$, $g : \RR
  \to \RR_+$ a Lipschitz function and $\tmmathbf{a}: \Ssurf
  \to \RR^3$ a Lipschitz vector field.
  
  \medskip \noindent If $ \mathfrak{h}_1 (t) \JW (t) > \sqrt{2 \pi}$ for every $t\in I$, the following assertions hold:
  \begin{enumerateroman}
    \item If the vector field $\tmmathbf{a}$ is axially symmetric, then any
    minimizer $\jm \in H^1 ( \Ssurf, \Tsurf )$ of
    $\mathcal{E}_{\jom}$ has an axially symmetric representative, in the sense
    that if $\jm$ is a minimizer of $\mathcal{E}_{\jom}$ then there exists an
    axially symmetric vector field $\ju \in H^1 ( \Ssurf, \Tsurf
    )$, built from $\jm$, such that $\mathcal{E}_{\jom} ( \ju
    ) =\mathcal{E}_{\jom} ( \jm )$.
    
    \item If the vector field $\tmmathbf{a}$ is axially antisymmetric, then
    any minimizer $\jm \in H^1 ( \Ssurf, \Tsurf )$ of
    $\mathcal{E}_{\jom}$ has an axially antisymmetric representative, in the
    sense that if $\jm$ is a minimizer of $\mathcal{E}_{\jom}$ then there
    exists an axially symmetric vector field $\ju \in H^1 ( \Ssurf,
    \Tsurf)$, built from $\jm$, such that $\mathcal{E}_{\jom} (
    \ju ) =\mathcal{E}_{\jom} ( \jm )$.
  \end{enumerateroman}
  In any case, if $\tmmathbf{a}$ is axially symmetric or antisymmetric,
  when $\mathfrak{h}_1 (t) \JW (t) > \sqrt{2 \pi}$ any minimizer $\jm \in H^1
  ( \Ssurf, \Tsurf )$ of $\mathcal{E}_{\jom}$ is of the form
  \begin{equation}
    \jm (\phi, t) =\tmmathbf{\alpha}_{\bot} (t) \cos \phi
    +\tmmathbf{\beta}_{\bot} (t) \sin \phi + \eta (t) \ee_3 .
    \label{eq:genexprmperp}
  \end{equation}
  for some $\tmmathbf{\alpha}_{\bot}, \tmmathbf{\beta}_{\bot} \in H^1(
  I, \RR^2)$, $\eta \in H^1( I, \RR )$, subject to the
  constraint of the resulting $\jm$ being $\Tsurf$-valued.
\end{theorem}

\begin{remark}
  Configurations of the form \eqref{eq:genexprmperp} belong to the class of vector fields satisfying the property that 
$\langle \jm_{\bot}
  \rangle_{\Ssurf_{\xi}} = 0$ for every $\xi \in \Ssurf$, and later on referred to as \emph{axially null-average} (see Definition \ref{def:axiallynullav}). Thus, in particular, Theorem~\ref{thm:main0} says that any minimizer is axially null-average  provided that $\mathfrak{h}_1 (t) \JW (t) > \sqrt{2 \pi}$ for every $t\in I$.
    When $\inf_{t\in I} (\mathfrak{h}_1 (t) \JW (t)) = \sqrt{2 \pi}$ we cannot conclude  that minimizers are still axially
  null-average. However, it will be evident from the proof
  ({\tmabbr{cf.}}~\eqref{eq:dphizerozero}) that any minimizer must be of the (more general) form $\jm (\phi, t)
  =\tmmathbf{\zeta} (t) +\tmmathbf{\alpha}_{\bot} (t) \cos \phi
  +\tmmathbf{\beta}_{\bot} (t) \sin \phi$ for suitable
  $\tmmathbf{\alpha}_{\bot}, \tmmathbf{\beta}_{\bot} \in H^1 ( I, \RR^2
  )$, $\tmmathbf{\zeta} \in H^1 ( I, \RR^3 )$ subject to the
  constraint of the resulting $\jm$ being $\Tsurf$-valued.
\end{remark}

\begin{remark}
  The assumption of $\Tsurf$ being a surface of revolution (around the
  $\ee_3$-axis) cannot be removed. Indeed, the entire argument is based on the
  closure property that if $\jm$ is $\Tsurf$-valued, then $\Rm^{\top} (\phi)
  \jm (\xi) \in \Tsurf$ for every $\phi \in \RR$ and every $\xi \in \Ssurf$.
  However, note that it is not necessary to (and we do not) assume
  $\tmmathbf{a}$ to be $\Tsurf$-valued. When we say that $\tmmathbf{a}: \Ssurf
  \to \RR^3$ is, e.g., axially symmetric, what we rigorously mean,
  that matches with Definition~\ref{def:axialsymmetry}, is that the image of
  $\tmmathbf{a}$ is included in a surface of revolution which can be different
  from $\Tsurf$.
\end{remark}

\begin{remark}
  As it will be apparent from the proofs, our arguments still work if one
  considers a multivariable potential of the form $g (\tmmathbf{m} (\xi) \cdot
  \tmmathbf{a}_1 (\xi), \tmmathbf{m} (\xi) \cdot \tmmathbf{a}_2 (\xi), \ldots,
  \tmmathbf{m} (\xi) \cdot \tmmathbf{a}_n (\xi))$ where the vector fields
  $(\tmmathbf{a}_i)_{i = 1}^n$ are either all axially symmetric or all axially
  antisymmetric. Also, as we explain in Remark~\ref{rmk:boundaryconditions},
  our arguments extend to the case in which the surface of revolution $\Ssurf$
  has a boundary, and we look for minimizers in $H^1 ( \Ssurf, \Tsurf
  )$ satisfying prescribed axially symmetric (or axially antisymmetric)
  Dirichlet boundary conditions. More generally, our results still hold when
  the anisotropy energy $\mathcal{A}_{\Ssurf}$ assumes the form
  \[ \mathcal{A}_{\Ssurf} ( \jm ) \assign \int_{\Ssurf} g
     (\jm (\xi) \cdot \tmmathbf{a}_1 (\xi)) \mathd \xi +
     \int_{\partial \Ssurf} b ( \jm (\sigma) \cdot \tmmathbf{a}_2
     (\sigma)) \mathd \sigma, \]
  i.e., when, as in {\cite{Di_Fratta_2022}}, boundary penalizations are
  involved. To control the complexity of the formulas, we do not treat these
  generalizations. However, it should be evident from the arguments that our
  results cover these more general settings with no technical changes at all.
\end{remark}

The expression of the minimizer $\jm$ in \eqref{eq:genexprmperp} allows for
both axially symmetric and antisymmetric vector fields. Indeed, if we set
$\tmmathbf{\alpha} (t) \assign \tmmathbf{\alpha}_{\bot} (t) + \eta (t) \ee_3$
then we can express the vector field $\jm$ in \eqref{eq:genexprmperp} under
the equivalent forms
({\tmabbr{cf.}}~\eqref{eq:pertaxsym}-\eqref{eq:pertaxantisym})
\begin{equation}
  \jm (\phi, t) \eqs A^{\T} (\phi) \tmmathbf{\alpha} (t) + (\sin \phi) (
  \tmmathbf{\beta}_{\bot} (t) - \ee_3 \times \tmmathbf{\alpha}_{\bot} (t)
  ),
\end{equation}
and
\begin{equation}
  \jm (\phi, t) \eqs A (\phi) \tmmathbf{\alpha} (t) + (\sin \phi) (
  \tmmathbf{\beta}_{\bot} (t) + \ee_3 \times \tmmathbf{\alpha}_{\bot} (t)
  ) .
\end{equation}
Thus, $\jm$ is axially symmetric when $\tmmathbf{\beta}_{\bot} (t) = \ee_3
\times \tmmathbf{\alpha}_{\bot} (t)$, and it is axially antisymmetric when
$\tmmathbf{\beta}_{\bot} (t) = - \ee_3 \times \tmmathbf{\alpha}_{\bot} (t)$.
Theorem~\ref{thm:main0} guarantees the existence of axially symmetric or
axially antisymmetric representatives, but, in general, it cannot rule out the
coexistence of minimizers with a broken symmetry term of the type $(\sin \phi)
( \tmmathbf{\beta}_{\bot} (t) \pm \ee_3 \times \tmmathbf{\alpha}_{\bot}
(t) )$. But things go better if one requires an additional condition on
the target manifold $\Tsurf$, which prevents the existence of flat zones in
$\Tsurf$.

\begin{definition}
  \label{def:neverflat}Let $\tmmathbf{\pi}_{\bot}$ be the projection of
  $\RR^3$ onto $\Sigma_0 \assign \RR^2 \times \{ 0 \}$ and, for every $z \in
  \RR$, \ let $\Sigma_z \assign z \ee_3 + \Sigma_0$ be the plane parallel to
  $\Sigma_0$ and passing through $z \ee_3$. Note that if $\Tsurf$ is a surface
  of revolution, then $\tmmathbf{\pi}_{\bot}( \Tsurf \cap \Sigma_z
  )$ is either empty or the union of circles in $\Sigma_0$. We say that
  $\Tsurf$ is {\tmem{never flat}} if for any $z \in \RR$ the following
  property holds: either $\tmmathbf{\pi}_{\bot} ( \Tsurf \cap \Sigma_z
  )$ is empty, or $\tmmathbf{\pi}_{\bot} ( \Tsurf \cap \Sigma_z
  )$ consists of a unique circle, or it consists of a finite family of
  circles at a positive distance from each other.
\end{definition}

If the profile $\jgam$ generating the surface of revolution $\Tsurf$ is the graph of a smooth function, i.e., of the form $\jgam(t)=(t,0,z(t))$, then the never-flat condition amounts to asking that the function $z$ does not have intervals where it is constant.
%
%
The never-flat condition is quite general and satisfied, e.g., by all
surfaces of revolution represented in Figure~\ref{fig:surfsrevol}. Instead,
the unit disk in $\Sigma_0$ does not satisfy the never-flat condition.

\begin{theorem}
  \label{thm:main1}For given surfaces of revolution $\Ssurf, \Tsurf \subseteq
  \RR^3$, we consider the energy functional defined for every $\jm \in H^1
  ( \Ssurf, \Tsurf )$ by
  {\opt}cf.{\tmem{~\eqref{eq:mainenfunc}}}{\cpt}
  \[ \mathcal{E}_{\jom} ( \jm ) \assign \int_{\Ssurf} | \grad
     \jm (\xi) |^2 \mathd \xi + \int_{\Ssurf} g ( \jm (\xi) \cdot
     \tmmathbf{a} (\xi) ) \mathd \xi + \int_{\Ssurf} | \langle
     \jm \rangle_{\Ssurf_{\xi}} \times \ee_3|^2  \jom^2 (\xi)
     \mathd \xi \]
  with $\jom : \Ssurf \to \RR_+$ a measurable weight such that $\sup_{t\in I} 
  (\mathfrak{h}_1 (t) \JW (t))<+\infty$, $g : \RR
  \to \RR_+$ a Lipschitz function and $\tmmathbf{a}: \Ssurf
  \to \RR^3$ a Lipschitz vector field.
  
   \medskip \noindent If $ \mathfrak{h}_1 (t) \JW (t) > \sqrt{2 \pi}$ for every $t\in I$ and $\Tsurf$ is never flat, then the following
  assertions hold:
  \begin{enumerateroman}
    \item If the vector field $\tmmathbf{a}$ is axially symmetric, then either
    every minimizer is axially symmetric, or there also coexist minimizers
    with line symmetry. An axially symmetric minimizer is given by $\jm (\phi,
    t) \assign A^{\T} (\phi) \tmmathbf{\gamma}_s (t)$ where
    $\tmmathbf{\gamma}_s$ \ is a solution to the one-dimensional minimization
    problem
    \[ \mathcal{F}_s (\tmmathbf{\gamma}) \assign \mathcal{E}_{\jom} (
       A^{\T} (\phi) \tmmathbf{\gamma} (t) ) \]
    among all possible profiles $\tmmathbf{\gamma} \in H^1 ( I, \Tsurf
    )$.
    
    \item If the vector field $\tmmathbf{a}$ is axially antisymmetric, then
    either every minimizer is axially antisymmetric, or there also coexist
    minimizers with line symmetry. An axially antisymmetric minimizer is given
    by $\jm (\phi, t) \assign A (\phi) \tmmathbf{\gamma}_a (t)$ where
    $\tmmathbf{\gamma}_a$ \ is a solution to the one-dimensional minimization
    problem
    \[ \mathcal{F}_s (\tmmathbf{\gamma}) \assign \mathcal{E}_{\jom} (A (\phi)
       \tmmathbf{\gamma} (t)) \]
    among all possible profiles $\tmmathbf{\gamma} \in H^1 ( I, \Tsurf
    )$.
  \end{enumerateroman}
  In any case, if the vector field $\tmmathbf{a}$ is axially symmetric or
  antisymmetric, then any minimizer $\jm \in H^1( \Ssurf, \Tsurf
  )$ of $\mathcal{E}_{\jom}$ has line symmetry and is of the form
  \begin{equation}
    \jm (\phi, t) =\tmmathbf{\alpha}_{\bot} (t) \cos \phi
    +\tmmathbf{\beta}_{\bot} (t) \sin \phi + \eta (t) \ee_3 .
    \label{eq:genexprmperpline}
  \end{equation}
  for some $\tmmathbf{\alpha}_{\bot}, \tmmathbf{\beta}_{\bot} \in H^1 (
  I, \RR^2 )$, $\eta \in H^1 ( I, \RR)$, subject to the
  orthogonality conditions $\tmmathbf{\beta}_{\bot} (t) \cdot
  \tmmathbf{\alpha}_{\bot} (t) \equiv 0$ and $| \tmmathbf{\beta}_{\bot} (t) |
  \equiv | \tmmathbf{\alpha}_{\bot} (t) |$, and to the conditions of the
  resulting $\jm$ being $\Tsurf$-valued.
\end{theorem}
\section{Axially (anti)symmetric minimizers: Proofs of Theorems
\ref{thm:main0} and \ref{thm:main1}}\label{sec:proofs}
\noindent We prove Theorems \ref{thm:main0} and \ref{thm:main1}, assuming that the
vector field $\tmmathbf{a}: \Ssurf \to \RR^3$ is axially symmetric.
The proof works the same when $\tmmathbf{a}: \Ssurf \to \RR^3$ is
axially antisymmetric, provided that the words {\tmem{symmetric}} and
{\tmem{antisymmetric}} are exchanged in the proper obvious places and formula
\eqref{eq:newmagtext} is replaced with the assignment $\ju (\phi, t) \assign
\Rm (\phi) \Rm^{\T} (\phi_{\ast}) \jm (\phi_{\ast}, t)$.

For clarity, we subdivide the proof in several steps. In order to improve the
readability of the argument, we summarize here the steps. In
{\tmname{Step~1}}, we reformulate the energy functional in local coordinates
through the parameterization $\xip$ of $\Ssurf$ defined by \eqref{eq:paramS}.
In {\tmname{Step~2}}, we derive an auxiliary estimate which, in particular,
assures that the construction in {\tmname{Step~3}} is well-posed. In
{\tmname{Step~3}}, we show that if $\mathfrak{h}_1 (t)\JW (t)
\geqslant \sqrt{2 \pi}$, then given any minimizer $\jm$ of
$\mathcal{E}_{\jom}$, there exists an axially symmetric vector field $\ju \in
H^1 ( \Ssurf, \Tsurf )$, built from $\jm$, with the same minimal
energy. In working out the details in {\tmname{Step~3}} we assume that the
minimal vector field $\jm$ is smooth, and we show in {\tmname{Step~5}} how to
avoid this assumption. In {\tmname{Step~4}}, we finalize the proofs of
Theorems \ref{thm:main0} and \ref{thm:main1}.

\medskip

{\noindent}{{\tmname{\begin{tabular}{|l|}
  \hline
  Step~1\\
  \hline
\end{tabular}}}
First, we rephrase the problem in local coordinates by
using the parameterization $\xip$ of $\Ssurf$ defined in \eqref{eq:paramS}.
With the usual abuse of notation, we set $\jm (\phi, t) \assign \jm (
\xip (\phi, t) )$, being aware that the context will always clarify such
an overload of the symbol $\jm$. In local coordinates, the shape anisotropy
term $\mathcal{A}_{\Ssurf}$ in \eqref{eq:energyterms0} reads as
\begin{equation}
  \mathcal{A}_{\Ssurf} ( \jm )  \eqs  \int_{t \in I} \int_0^{2
  \pi} g ( \jm (\phi, t) \cdot \tmmathbf{a} (\phi, t) ) \mathd \phi
  \sqrt{\mathfrak{g} (t)} \mathd t. 
\end{equation}
Also, with $\jm_{\bot} \assign ( \jm \cdot \ee_1 ) \ee_1 + (
\jm \cdot \ee_2 ) \ee_2$, and the assignment $\langle \jm_{\bot}
\rangle (t) \assign \frac{1}{2 \pi} \int_0^{2 \pi} \jm_{\bot} (\theta,
t) \mathd \theta$, the  term $\mathcal{P}_{\Ssurf}$ in
\eqref{eq:energyterms1} can be written as follows:
\begin{align}
  \mathcal{P}_{\Ssurf} ( \jm ) & =  \int_{\Ssurf} | \jom
  (\xi) \langle \jm_{\bot} \rangle_{\Ssurf_{\xi}} |^2 \mathd
  \xi \\
  & =  \int_{t \in I} \int_0^{2 \pi} \left| \, \frac{\jom ( \xip (\phi,
  t) )}{| \Ssurf_{\xip (\phi, t)} |} \int_{\Ssurf_{\xip
  (\phi, t)}} \jm_{\bot} (\sigma) \mathd \sigma \, \right|^2 \mathd \phi
  \sqrt{\mathfrak{g} (t)} \mathd t \\
  & \eqs  \int_{t \in I} \int_0^{2 \pi} \jom^2 ( \xip (\phi, t) )
  \left| \, \frac{1}{| 2 \pi \mathfrak{h}_1(t) |} \int_0^{2 \pi} \jm_{\bot}
  (\theta, t)  | \partial_{\theta} \xip (\theta, t) | \mathd \theta
  \,\right|^2 \mathd \phi \sqrt{\mathfrak{g} (t)} \mathd t \\
  & \eqs  \int_{t \in I} \left( \int_0^{2 \pi} \jom^2 ( \xip (\phi, t)
  ) \mathd \phi \right)  | \langle \jm_{\bot} \rangle
  (t) |^2  \sqrt{\mathfrak{g} (t)} \mathd t \\
  & =  \int_{t \in I} \JW^2 (t)  | \langle \jm_{\bot}
  \rangle (t) |^2  \sqrt{\mathfrak{g} (t)} \mathd t. 
\end{align}
In computing the previous expressions we took into account that $| \Ssurf_{\xip (\phi, t)} |=2 \pi \mathfrak{h}_1(t)$ and $| \partial_{\theta} \xip (\theta, t) |= \mathfrak{h}_1(t)$.

Finally, we focus on the expression in local coordinates of the Dirichlet
energy term $\mathcal{D}_{\Ssurf}$ in \eqref{eq:energyterms0}. To that end, we
observe that with the notation introduced in Section~\ref{sec:setup}, in
particular \eqref{eq:tauphi} and \eqref{eq:taut}, and with the usual abuse of notation
$\jm (\phi, t) \assign ( \jm \circ \xip ) (\phi, t)$,
 there holds that
\begin{equation}
  ( \partial_{\tfi} \jm \circ \xip ) (\phi, t)  \eqs
  \partial_{\phi} \jm (\phi, t), \qquad ( \partial_{\tti} \jm \circ \xip
  ) (\phi, t)  \eqs \partial_t \jm (\phi, t),
\end{equation}
because, e.g., $\partial_{\phi} \jm (\phi, t) = [ D \jm \circ \xip (\phi,
t) ] \tfi (\phi, t)$. Hence, we find that{\tmsamp{}}
\begin{equation}
  | ( \nabla \jm \circ \xip ) (\phi, t) |^2   \eqs 
  \frac{| \partial_{\phi} \jm (\phi, t) |^2}{\mathfrak{h}_1^2 (t)}
  + \frac{| \partial_t \jm (\phi, t) |^2}{\mathfrak{h}_2^2 (t)} 
\end{equation}
with $\mathfrak{h}_1$ and $\mathfrak{h}_2$ the metric coefficients defined in
\eqref{eq:metrich1}. It follows that the Dirichlet energy reads under the form
\begin{equation}
  \int_{\Ssurf} | \nabla \jm (\xi) |^2 \mathd \xi  \eqs  \int_{t
  \in I} \int_0^{2 \pi} \frac{| \partial_{\phi} \jm (\phi, t)
  |^2}{\mathfrak{h}^2_1 (t)} + \frac{| \partial_t \jm (\phi, t)
  |^2}{\mathfrak{h}_2^2 (t)} \mathd \phi \sqrt{\mathfrak{g} (t)} \mathd
  t. 
\end{equation}
Overall, in the coordinate chart induced by $\xip$, the energy functional on
$H^1 ( \Ssurf, \Tsurf )$ assumes the form
\begin{align}
  \mathcal{E}_{\jom} ( \jm ) & \eqs  \int_{t \in I} \int_0^{2 \pi}
  \frac{| \partial_{\phi} \jm (\phi, t) |^2}{\mathfrak{h}^2_1 (t)}
  + \frac{| \partial_t \jm (\phi, t) |^2}{\mathfrak{h}^2_2 (t)}
  \mathd \phi \sqrt{\mathfrak{g} (t)} \mathd t \nonumber\\
  &   \quad \quad \quad \quad + \int_{t \in I} \JW^2 (t)  |
  \langle \jm_{\bot} \rangle (t) |^2  \sqrt{\mathfrak{g} (t)}
  \mathd t \nonumber\\
  &   \quad \quad \quad \quad \quad \quad \quad \quad \quad + \int_{t \in I}
  \int_0^{2 \pi} g ( \jm (\phi, t) \cdot \tmmathbf{a} (\phi, t) )
  \mathd \phi \sqrt{\mathfrak{g} (t)} \mathd t.  \label{eq:Einloccoords}
\end{align}
From now on, we are going to work with the local coordinates expression \eqref{eq:Einloccoords}.

\medskip

{\noindent}{{\tmname{\begin{tabular}{|l|}
  \hline
  Step~2\\
  \hline
\end{tabular}}}
The main aim of this step is to guarantee that the
construction in the next {\tmname{Step~3}} is well-posed. In what follows, we denoted by
$\jx (\xi)$ the distance of $\xi \in \Ssurf$ from the $\ee_3$-axis, i.e., the
radius of the circle of latitude $\Ssurf_{\xi}$. Note that, in local
coordinates, $\jx ( \xip (\phi, t) ) =\mathfrak{h}_1 (t)$. We want
to show that if $\jm \in H^1 ( \Ssurf, \Tsurf )$ has finite
$\mathcal{E}_{\jom}$-energy (in particular, if it is a minimizer) and
$\mathfrak{h}^2_1 (t) \JW^2 (t) \geqslant 2 \pi$ then
\begin{equation}
  \int_{\Ssurf} | \jm_{\bot} (\xi) |^2  \frac{\mathd \xi}{\jx^2
  (\xi)} < + \infty . \label{eq:xmperpfinite}
\end{equation}
In order to establish \eqref{eq:xmperpfinite}, it is sufficient to show that
in local coordinates one has
\begin{equation}
  \int_{t \in I} \int_0^{2 \pi} \frac{| \jm_{\bot} (\phi, t)
  |^2}{\mathfrak{h}^2_1 (t)} \mathd \phi \sqrt{\mathfrak{g} (t)} \mathd
  t < + \infty \label{eq:xmperpfiniteloccoord} .
\end{equation}
For that, we recall that for periodic functions there holds the following
Poincar{\'e}-Wirtinger inequality
\begin{equation}
  \int_0^{2 \pi} | \ju (\phi) - \langle \ju \rangle |^2
  \mathd \phi \leqslant \int_0^{2 \pi} | \partial_{\phi} \ju (\phi)
  |^2 \mathd \phi, \quad \langle \ju \rangle \assign (2
  \pi)^{^{- 1}} \int_0^{2 \pi} \ju (\phi) \mathd \phi . \label{eq:PWforper}
\end{equation}
Moreover, given the quadratic setting, the left-hand side of the previous
relation can be written as
\begin{equation}
  \int_0^{2 \pi} | \ju (\phi) - \langle \ju \rangle |^2
  \mathd \phi \eqs \int_0^{2 \pi} | \ju (\phi) |^2 - |
  \langle \ju \rangle |^2 \mathd \phi .
\end{equation}
Therefore, if we write $\jm = \jm_{\bot} + ( \jm \cdot \ee_3 )
\ee_3$ then for every $t \in I$, there holds
\begin{align}
  0 & \leqslant  \int_{t \in I} \int_0^{2 \pi} \frac{| \jm_{\bot} (\phi,
  t) |^2 - | \langle \jm_{\bot} \rangle (t)
  |^2}{\mathfrak{h}^2_1 (t)} \mathd \phi \sqrt{\mathfrak{g} (t)} \mathd
  t \nonumber\\
  & \leqslant  \int_{t \in I} \int_0^{2 \pi}  \frac{| \partial_{\phi}
  \jm_{\bot} (\phi, t) |^2}{\mathfrak{h}^2_1 (t)} \mathd \phi
  \sqrt{\mathfrak{g} (t)} \mathd t. 
\end{align}
Hence, whenever $\mathfrak{h}^2_1 (t) \JW^2 (t) \geqslant 2 \pi$ we get
\begin{align}
  0 & \leqslant  \int_{t \in I} \int_0^{2 \pi} \frac{| \jm_{\bot} (\phi,
  t) |^2}{\mathfrak{h}^2_1 (t)} \mathd \phi \sqrt{\mathfrak{g} (t)}
  \mathd t  \label{eq:mperpfinite}\\
  & \leqslant  \int_{t \in I} \int_0^{2 \pi}  \frac{| \partial_{\phi}
  \jm_{\bot} (\phi, t) |^2}{\mathfrak{h}^2_1 (t)} \mathd \phi
  \sqrt{\mathfrak{g} (t)} \mathd t + 2 \pi \int_{t \in I} \frac{|
  \langle \jm_{\bot} \rangle (t) |^2}{\mathfrak{h}^2_1 (t)}
  \sqrt{\mathfrak{g} (t)} \mathd t \\
  & \leqslant  \int_{t \in I} \int_0^{2 \pi}  \frac{| \partial_{\phi}
  \jm_{\bot} (\phi, t) |^2}{\mathfrak{h}^2_1 (t)} \mathd \phi
  \sqrt{\mathfrak{g} (t)} \mathd t + \int_{t \in I} \JW^2 (t) |
  \langle \jm_{\bot} \rangle (t) |^2 \sqrt{\mathfrak{g} (t)}
  \mathd t,  \label{eq:mperpfinite2}
\end{align}
and this shows that the integral in \eqref{eq:xmperpfinite} is finite provided
that $\jm \in H^1 ( \Ssurf, \Tsurf )$ has finite
$\mathcal{E}_{\jom}$-energy and $\mathfrak{h}^2_1 (t) \JW^2 (t) \geqslant 2
\pi$.

\medskip

{\noindent}{{\tmname{\begin{tabular}{|l|}
  \hline
  Step~3\\
  \hline
\end{tabular}}} {\tmem{Given any minimizer $\jm$ of $\mathcal{E}_{\jom}$,
there exists an axially symmetric vector field $\ju \in H^1 ( \Ssurf,
\Tsurf )$, built from $\jm$, with the same minimal energy.}}

By direct methods in the calculus of variations, we know that there exists a (global) minimizer $\jm \in H^1 ( \Ssurf, \Tsurf )$ of the
energy functional $\mathcal{E}_{\jom}$, but, in general, more than one. We
want to show that given any (global) minimizer $\jm$ of $\mathcal{E}_{\jom}$,
there exists an axially symmetric texture $\ju \in H^1 ( \Ssurf, \Tsurf
)$, built from $\jm$, with the same minimal energy, i.e., such that
$\mathcal{E}_{\jom}( \ju ) =\mathcal{E}_{\jom} ( \jm
)$. 
We
first assume that $\jm \in H^1 ( \Ssurf, \Tsurf )$ is a smooth
minimizer of the energy $\mathcal{E}_{\jom}$. Afterward, in {\tmname{Step~5}},
we show how to remove this assumption through a density argument.

The idea is to introduce the real variable function
\begin{equation}
  \Phi_{\mathcal{E}} : \phi \in [0, 2 \pi] \mapsto \; {\Phi_{\mathcal{D}}} 
  (\phi) + \Phi_{\mathcal{A}} (\phi) \in \RR_+, \label{eq:Phifull}
\end{equation}
where
\begin{align}
  {\Phi_{\mathcal{D}}}  (\phi) & \assign  \int_{t \in I} \left( \frac{|
  \jm_{\bot} (\phi, t) |^2}{\mathfrak{h}^2_1 (t)} + \frac{|
  \partial_t \jm (\phi, t) |^2}{\mathfrak{h}^2_2 (t)}  \right)
  \sqrt{\mathfrak{g} (t)} \mathd t, \\
  \Phi_{\mathcal{A}} (\phi) & \assign  \int_{t \in I} g ( \jm (\phi, t)
  \cdot \tmmathbf{a} (\phi, t) )  \sqrt{\mathfrak{g} (t)} \mathd t . 
\end{align}
For $\phi_{\ast} \in \tmop{argmin}_{\phi \in [0, 2 \pi]} \Phi_{\mathcal{E}}
(\phi)$ we have that
\begin{align}
  \Phi_{\mathcal{E}} (\phi_{\ast}) & \eqs  \int_{t \in I} \left( \frac{|
  \jm_{\bot} (\phi_{\ast}, t) |^2}{\mathfrak{h}^2_1 (t)} + \frac{|
  \partial_t \jm (\phi_{\ast}, t) |^2}{\mathfrak{h}^2_2 (t)} \right) 
  \sqrt{\mathfrak{g} (t)} \mathd t \nonumber\\
  &   \qquad \qquad \qquad \qquad + \int_{t \in I} g ( \jm
  (\phi_{\ast}, t) \cdot \tmmathbf{a} (\phi_{\ast}, t) )
  \sqrt{\mathfrak{g} (t)} \mathd t.  \label{eq:armingphifull}
\end{align}
Next, we define a new vector field $\ju \in H^1 ( \Ssurf, \Tsurf )$
via the relation
\begin{equation}
  \ju (\phi, t) \assign \Rm^{\top} (\phi) \Rm (\phi_{\ast}) \jm (\phi_{\ast},
  t) . \label{eq:newmagtext}
\end{equation}
We then have $ \partial_t \ju (\phi, t)  \eqs  \Rm^{\top} (\phi) \Rm (\phi_{\ast})
  \partial_t \jm (\phi_{\ast}, t)$ and $\partial_{\phi} \ju (\phi, t)  \eqs  \partial_{\phi} \Rm^{\top} (\phi) \Rm
  (\phi_{\ast}) \jm (\phi_{\ast}, t)$,
from which, observing that $\Rm (\phi) \partial_{\phi} \Rm^{\T} (\phi)
\tmmathbf{v}= \ee_3 \times \tmmathbf{v}$ for every $\tmmathbf{v} \in \RR^3$,
we infer that
\begin{equation}
  | \partial_t \ju (\phi, t) |^2  \eqs  | \partial_t \jm
  (\phi_{\ast}, t) |^2, 
\end{equation}
and
\begin{align}
  | \partial_{\phi} \ju (\phi, t) |^2 & \eqs  | \Rm (\phi)
  \partial_{\phi} \Rm^{\top} (\phi) \Rm (\phi_{\ast}) \jm (\phi_{\ast}, t)
  |^2 \nonumber\\
  & \eqs  | \ee_3 \times \Rm (\phi_{\ast}) \jm (\phi_{\ast}, t)
  |^2 \nonumber\\
  & \eqs  | \ee_3 \times \jm (\phi_{\ast}, t) |^2 \nonumber\\
  & \eqs  | \jm_{\bot} (\phi_{\ast}, t) |^2 . 
\end{align}
\begin{remark}
  ({\tmem{The construction preserves boundary conditions}})
  \label{rmk:boundaryconditions}The construction that leads to
  \eqref{eq:newmagtext} is compatible with eventually prescribed boundary
  conditions: if $\Ssurf$ is a surface of revolution with boundary, then
  $\jm_{\partial \Ssurf} \equiv \ju_{\partial \Ssurf}$. Indeed, suppose to fix
  the ideas that $I = [0, 1]$ and that, e.g., $\jm (\phi, 1) = \Rm^{\top}
  (\phi) \ee$ for some $\ee \in \Tsurf$, then we have $\ju (\phi, 1) \assign
  \Rm^{\top} (\phi) \Rm (\phi_{\ast}) \jm (\phi_{\ast}, 1) = \Rm^{\top} (\phi)
  \Rm (\phi_{\ast}) \Rm^{\top} (\phi_{\ast}) \ee = \Rm^{\top} (\phi) \ee$ and,
  therefore, the boundary value is preserved.
\end{remark}
Also, using the local coordinate representation of the characterization of
axially symmetric vector fields expressed by
$\eqref{eq:coordinatefreeaxsymm}$, we obtain that
\begin{align}
  g ( \ju (\phi, t) \cdot \tmmathbf{a} (\phi, t) ) & \eqs  g
  ( \Rm^{\top} (\phi) \Rm (\phi_{\ast}) \jm (\phi_{\ast}, t) \cdot
  \Rm^{\top} (\phi) \tmmathbf{a} (0, t) ) \nonumber\\
  & \eqs  g ( \Rm (\phi_{\ast}) \jm (\phi_{\ast}, t) \cdot \tmmathbf{a}
  (0, t) ) \nonumber\\
  & \eqs  g ( \jm (\phi_{\ast}, t) \cdot \Rm^{\top} (\phi_{\ast})
  \tmmathbf{a} (0, t) ) \nonumber\\
  & \eqs  g ( \jm (\phi_{\ast}, t) \cdot \tmmathbf{a} (\phi_{\ast}, t)
  ) . 
\end{align}
From the previous computations, using $\eqref{eq:armingphifull}$ and
$\eqref{eq:mperpfinite2}$, we infer that when $\mathfrak{h}^2_1 (t)
\JW^2 (t) \geqslant 2 \pi$, the following estimates hold
\begin{align}
  \mathcal{E}_{\jom}( \ju ) & \eqs  \int_0^{2 \pi} \int_{t \in I}
  \frac{| \jm_{\bot} (\phi_{\ast}, t) |^2}{\mathfrak{h}^2_1 (t)} +
  \frac{| \partial_t \jm (\phi_{\ast}, t) |^2}{\mathfrak{h}^2_2
  (t)}  \sqrt{\mathfrak{g} (t)} \mathd t \mathd \phi \nonumber\\
  &   \qquad \qquad + \int_0^{2 \pi} \int_{t \in I} g (
  \jm (\phi_{\ast}, t) \cdot \tmmathbf{a} (\phi_{\ast}, t) ) 
  \sqrt{\mathfrak{g} (t)} \mathd t \mathd \phi \\
  & \eqs  \int_0^{2 \pi}
  \Phi_{\mathcal{E}} (\phi_{\ast}) \mathd \phi \\
  & \leqslant  \int_0^{2 \pi} \Phi_{\mathcal{E}} (\phi) \mathd \phi \\
  & \eqs  \int_{t \in I} \int_0^{2 \pi} \frac{| \jm_{\bot} (\phi, t)
  |^2}{\mathfrak{h}^2_1 (t)} + \frac{| \partial_t \jm (\phi, t)
  |^2}{\mathfrak{h}^2_2 (t)} \sqrt{\mathfrak{g} (t)} \mathd \phi \mathd
  t \nonumber\\
  &  \qquad \qquad  + \int_{t \in I} \int_0^{2 \pi} g (
  \jm (\phi, t) \cdot \tmmathbf{a} (\phi, t) ) \mathd \phi
  \sqrt{\mathfrak{g} (t)} \mathd t  \label{eq:equality1}\\
  & \leqslant  \int_{t \in I} \int_0^{2
  \pi}  \frac{| \partial_{\phi} \jm_{\bot} (\phi, t)
  |^2}{\mathfrak{h}^2_1 (t)} + \frac{| \partial_t \jm (\phi, t)
  |^2}{\mathfrak{h}^2_2 (t)} \mathd \phi \sqrt{\mathfrak{g} (t)} \mathd
  t \nonumber\\
  &   \qquad \qquad + \int_{t \in I} \JW^2 (t) |
  \langle \jm_{\bot} \rangle (t) |^2 \sqrt{\mathfrak{g} (t)}
  \mathd t \nonumber\\
  &   \qquad \qquad \qquad \qquad + \int_{t \in I} \int_0^{2
  \pi} g ( \jm (\phi, t) \cdot \tmmathbf{a} (\phi, t) ) \mathd \phi
  \sqrt{\mathfrak{g} (t)} \mathd t  \label{eq:equality2}\\
  & \leqslant  \mathcal{E}_{\jom} ( \jm ) .  \label{eq:equality3}
\end{align}
The minimality of $\jm$ entails the equality $\mathcal{E}_{\jom} ( \jm
) =\mathcal{E}_{\jom}( \ju )$ from which the conclusions of
Theorem~\ref{thm:main0} follow.

\medskip
{\noindent}{{\tmname{\begin{tabular}{|l|}
  \hline
  Step~4\\
  \hline
\end{tabular}}}
 Up to now, we know that if $\jm$ is a smooth minimizer of
$\mathcal{E}_{\jom}$ then there exists an axially symmetric vector field $\ju
\in H^1 ( \Ssurf, \Tsurf )$, built from $\jm$, such that
$\mathcal{E}_{\jom}( \ju ) =\mathcal{E}_{\jom} ( \jm
)$. We now show that, as a consequence, if $\mathfrak{h}_1 (t)
\JW (t) > \sqrt{2 \pi}$, then the conclusions of Theorem~\ref{thm:main1} hold,
i.e., that any minimizer of $\mathcal{E}_{\jom}$ is necessarily axially
symmetric. Indeed, given that $\mathcal{E}_{\jom}( \ju )
=\mathcal{E}_{\jom} ( \jm )$, one gets that
\eqref{eq:equality1}$\eqs$\eqref{eq:equality2}$=$\eqref{eq:equality3}, and
these equalities entail that
\begin{align}
  \int_{t \in I} \int_0^{2 \pi}  \frac{| \partial_{\phi} \jm (\phi, t)
  |^2}{\mathfrak{h}^2_1 (t)} & +  \int_{t \in I} \left[ \JW^2 (t) -
  \frac{2 \pi}{\mathfrak{h}^2_1 (t)} \right] | \langle \jm_{\bot}
  \rangle (t) |^2 \sqrt{\mathfrak{g} (t)} \mathd t \nonumber\\
  &   \qquad \eqs \int_{t \in I} \int_0^{2 \pi} \frac{| \jm_{\bot}
  (\phi, t) - \langle \jm_{\bot} \rangle (t)
  |^2}{\mathfrak{h}^2_1 (t)} \mathd \phi \sqrt{\mathfrak{g} (t)} \mathd
  t.  \label{eq:conseqminPW}
\end{align}
But then, Poincar{\'e}-Wirtinger inequality \eqref{eq:PWforper}, together with
the previous equality \eqref{eq:conseqminPW}, implies that any minimizer for
which \eqref{eq:equality1}$\eqs$\eqref{eq:equality2}$=$\eqref{eq:equality3},
must necessarily satisfy the relations
\begin{equation}
  \left[ \JW^2 (t) - \frac{2 \pi}{\mathfrak{h}^2_1 (t)} \right] |
  \langle \jm_{\bot} \rangle (t) |^2 = 0,
  \label{eq:dphizerozero}
\end{equation}
\begin{equation}
  | \partial_{\phi} \jm (\phi, t) \cdot \ee_3 |^2 \eqs 0.
  \label{eq:dphizero}
\end{equation}
In writing the previous two relations, we took into account that by hypotheses
(see Remark~\ref{rmk:regulsurf}), we have $\sqrt{\mathfrak{g} (t)} > 0$ for
every $t \in I$ with the possible exception of the two boundary points of
$\jgam$ if $\jgam$ touches the $\ee_3$-axis. Also, from
\eqref{eq:dphizerozero}, \eqref{eq:dphizero}, and \eqref{eq:conseqminPW} we
get that any minimizer of $\mathcal{E}_{{\jom}}$ satisfies the
relation
\begin{equation}
  \int_{t \in I} \left( \int_0^{2 \pi} ( | \partial_{\phi}
  \jm_{\bot} (\phi, t) |^2 - | \jm_{\bot} (\phi, t) - \langle
  \jm_{\bot} \rangle (t) |^2 ) \mathd \phi \right)
  \frac{\sqrt{\mathfrak{g} (t)}}{\mathfrak{h}^2_1 (t)} \mathd t \eqs 0,
\end{equation}
and the integrand is nonnegative by Poincar{\'e}--Wirtinger inequality. It
follows that for a.e. $t \in I$ there holds
\begin{equation}
  \int_0^{2 \pi} ( | \partial_{\phi} \jm_{\bot} (\phi, t) |^2
  - | \jm_{\bot} (\phi, t) - \langle \jm_{\bot} \rangle (t)
  |^2 ) \mathd \phi \eqs 0.
\end{equation}
But this means that the equality sign is reached in the Poincar{\'e}-Wirtinger
inequality \eqref{eq:PWforper}, and this is known to happen if, and only if,
$\jm_{\bot} (\phi, t) = \langle \jm_{\bot} \rangle (t)
+\tmmathbf{\alpha}_{\bot} (t) \cos \phi +\tmmathbf{\beta}_{\bot} (t) \sin
\phi$ for suitable functions $\tmmathbf{\alpha}_{\bot},
\tmmathbf{\beta}_{\bot} : t \in I \mapsto \RR^2 \times \{ 0 \}$. Also, we know
from \eqref{eq:dphizero} that $\jm (\phi, t) \cdot \ee_3$ depends only on the
$t$-variable. Therefore, every minimizer is of the form
\begin{equation}
  \jm (\phi, t) = \langle \jm_{\bot} \rangle (t)
  +\tmmathbf{\alpha}_{\bot} (t) \cos \phi +\tmmathbf{\beta}_{\bot} (t) \sin
  \phi + f (t) \ee_3 \label{eq:ansatzaposter}
\end{equation}
for a suitable scalar function $f : I \to \RR$ which is nothing but
$\jm (\phi, t) \cdot \ee_3$. Moreover, if $\mathfrak{h}_1 (t) \JW
(t) > \sqrt{2 \pi}$, then $\langle \jm_{\bot} \rangle (t) \equiv 0$
and, therefore, the minimizer $\jm$ is necessarily of the form
\begin{equation}
  \jm (\phi, t) =\tmmathbf{\alpha}_{\bot} (t) \cos \phi
  +\tmmathbf{\beta}_{\bot} (t) \sin \phi + f (t) \ee_3
  \label{eq:ansatzaposter2}
\end{equation}
Note that, in general, the previous expression includes both axially symmetric
and antisymmetric vector fields. Indeed, if we set $\tmmathbf{\alpha} (t)
\assign (\tmmathbf{\alpha}_{\bot} (t), f (t))$ then we can express the vector
field $\jm$ in \eqref{eq:ansatzaposter2} both as a perturbation of an axially
symmetric vector field,
\begin{align}
  \jm (\phi, t) & \eqs  (\cos \phi) \tmmathbf{\alpha}_{\bot} (t) + (\sin
  \phi) \ee_3 \times \tmmathbf{\alpha}_{\bot} (t) + ( \ee_3 \otimes \ee_3
  ) \tmmathbf{\alpha} (t) \nonumber\\
  &   \qquad \qquad \qquad + (\sin \phi) (
  \tmmathbf{\beta}_{\bot} (t) - \ee_3 \times \tmmathbf{\alpha}_{\bot} (t)
  ) \nonumber\\
  & \eqs  A^{\T} (\phi) \tmmathbf{\alpha} (t) + (\sin \phi) (
  \tmmathbf{\beta}_{\bot} (t) - \ee_3 \times \tmmathbf{\alpha}_{\bot} (t)
  )  \label{eq:pertaxsym}
\end{align}
and as a perturbation of an axially antisymmetric vector field,
\begin{align}
  \jm (\phi, t) & \eqs  (\cos \phi) \tmmathbf{\alpha}_{\bot} (t) - (\sin
  \phi) \ee_3 \times \tmmathbf{\alpha}_{\bot} (t) + ( \ee_3 \otimes \ee_3
  ) \tmmathbf{\alpha} (t) \nonumber\\
  &   \qquad \qquad \qquad + (\sin \phi) (
  \tmmathbf{\beta}_{\bot} (t) + \ee_3 \times \tmmathbf{\alpha}_{\bot} (t)
  ) \nonumber\\
  & \eqs  A (\phi) \tmmathbf{\alpha} (t) + (\sin \phi) (
  \tmmathbf{\beta}_{\bot} (t) + \ee_3 \times \tmmathbf{\alpha}_{\bot} (t)
  ) .  \label{eq:pertaxantisym}
\end{align}
In other words, $\jm$ is axially symmetric when $\tmmathbf{\beta}_{\bot} (t) =
\ee_3 \times \tmmathbf{\alpha}_{\bot} (t)$ and axially antisymmetric when
$\tmmathbf{\beta}_{\bot} (t) = - \ee_3 \times \tmmathbf{\alpha}_{\bot} (t)$.
As a consequence, for a minimal configuration $\jm$ of the form
\eqref{eq:ansatzaposter2} to have line symmetry, it is sufficient to satisfy
the orthogonality conditions
\begin{equation}
  | \tmmathbf{\beta}_{\bot} (t) | \equiv | \tmmathbf{\alpha}_{\bot} (t) |,
  \quad \tmmathbf{\beta}_{\bot} (t) \cdot \tmmathbf{\alpha}_{\bot} (t) \equiv
  0. \label{eq:orthnorm}
\end{equation}
Indeed, if these conditions are met, then for every $t \in I$ there holds
$\tmmathbf{\beta}_{\bot} (t) \eqs \pm \ee_3 \times \tmmathbf{\alpha}_{\bot}
(t)$, and therefore, from \eqref{eq:pertaxsym}-\eqref{eq:pertaxantisym}, that
for every $(\phi, t) \in \RR \times I$ either $\jm (\phi, t) = A^{\T} (\phi)
\tmmathbf{\alpha} (t)$ or $\jm (\phi, t) = A (\phi) \tmmathbf{\alpha} (t)$.

It remains to show that any minimizer satisfies the orthogonality conditions
\eqref{eq:orthnorm}. But this is a consequence of the assumption on the target
manifold $\Tsurf$ of being never flat
({\tmabbr{cf.}}~Definition~\ref{def:neverflat}). Indeed, for any $t \in I$ the
map $\jm_{\bot} (\cdot, t)$ takes values in $\tmmathbf{\pi}_{\bot} (
\Tsurf \cap \Sigma_{f (t)} )$. This means that $| \jm_{\bot}
(\cdot, t) |$ is a Sobolev function taking values into a finite set,
i.e., in the set of radii associated with the circles in
$\tmmathbf{\pi}_{\bot} ( \Tsurf \cap \Sigma_{f (t)})$. It follows
that $\phi \in I \mapsto | \jm_{\bot} (\cdot, t) |$ is constant
and, therefore, that $| \jm_{\bot} (\cdot, t) |$ depends only on
the $t$-variable. Therefore we get
\begin{equation}
  | \tmmathbf{\alpha}_{\bot} (t) |^2 = | \jm_{\bot} (0, t) |^2 =
  | \jm_{\bot} (\pi / 2, t) |^2 = | \tmmathbf{\beta}_{\bot} (t)
  |^2, \label{eq:firstorthog}
\end{equation}
and this proves the first condition in \eqref{eq:orthnorm}. But now, from
\eqref{eq:ansatzaposter2} and \eqref{eq:firstorthog}, we obtain
\begin{equation}
  | \jm_{\bot} (\phi, t) |^2 = | \tmmathbf{\alpha}_{\bot} (t) |^2
  +\tmmathbf{\beta}_{\bot} (t) \cdot \tmmathbf{\alpha}_{\bot} (t) \sin (2
  \phi)
\end{equation}
and the only way for $| \jm_{\bot} (\phi, t) |$ to be constant in
$\phi$ is that also the second orthogonality condition in \eqref{eq:orthnorm}
is satisfied.

If the vector field $\tmmathbf{a}$ is axially symmetric, given that we know
from {\tmname{Step~3}} the existence of axially symmetric minimizers of
$\mathcal{E}_{{\jom}}$, we conclude an axially symmetric
minimizer is given by $\jm (\phi, t) \assign A^{\T} (\phi) \tmmathbf{\gamma}_s
(t)$ where $\tmmathbf{\gamma}_s$ \ is a solution to the one-dimensional
minimization problem
\begin{align}
  \mathcal{F}_s (\tmmathbf{\gamma}) \assign \mathcal{E}_{\jom} ( A^{\T}
  (\phi) \tmmathbf{\gamma} (t) ) & \eqs  2 \pi \int_{t \in I} \frac{|
  \tmmathbf{\gamma}_{\bot} (t) |^2}{\mathfrak{h}^2_1 (t)} + \frac{|
  \dot{\tmmathbf{\gamma}}_{\bot} (t) |^2}{\mathfrak{h}^2_2 (t)} 
  \sqrt{\mathfrak{g} (t)} \mathd t \nonumber\\
  &   \qquad \qquad  + \int_{t \in I} \int_0^{2 \pi} g
  (\tmmathbf{\gamma} (t) \cdot A (\phi) \tmmathbf{a} (\phi, t)) \mathd \phi
  \sqrt{\mathfrak{g} (t)} \mathd t 
\end{align}
among all possible profiles $\tmmathbf{\gamma} \in H^1 ( I, \Tsurf
)$.

\medskip

{\noindent}{{\tmname{\begin{tabular}{|l|}
  \hline
  Step~5\\
  \hline
\end{tabular}}}
 In this final step, we show how to remove the smoothness
assumption on the minimizer. For that, we improve the ideas
in~{\cite{Fratta_2023}}. Let $\jm$ be a minimizer of
$\mathcal{E}_{{\jom}}$ and let $\jm^{\varepsilon} \in C^{\infty}
( \Ssurf, \Tsurf )$ be a family such that $\jm^{\varepsilon}
\to \jm$ strongly in $H^1 ( \Ssurf, \Tsurf )$. Such a
family certainly exists because of a well-known result of Schoen  and 
Uhlenbeck (see {\cite[p.267]{Schoen_1983}}). For every $\varepsilon > 0$, we
consider the continuous function
\begin{align}
  \Phi_{\mathcal{E}}^{\varepsilon} (\phi) & \eqs  \int_{t \in I} \left(
  \frac{| \jm^{\varepsilon}_{\bot} (\phi, t) |^2}{\mathfrak{h}^2_1
  (t)} + \frac{| \partial_t \jm^{\varepsilon} (\phi, t)
  |^2}{\mathfrak{h}^2_2 (t)} \right)  \sqrt{\mathfrak{g} (t)} \mathd t
  \nonumber\\
  &   \qquad \qquad \qquad \qquad + \int_{t \in I} g (
  \jm^{\varepsilon} (\phi, t) \cdot \tmmathbf{a} (\phi_{\ast}, t) )
  \sqrt{\mathfrak{g} (t)} \mathd t,  \label{eq:armingphifulleps}
\end{align}
and we choose a point $\phi_{\varepsilon} \in \tmop{argmin}_{\phi
\in [0, 2 \pi]} \Phi_{\mathcal{E}}^{\varepsilon} (\phi)$. We then define the
axially symmetric vector field $\ju_{\varepsilon} \in H^1 ( \Ssurf,
\Tsurf )$ via the relation
\begin{equation}
  \ju_{\varepsilon} (\phi, t) \assign \Rm^{\top} (\phi) \Rm
  (\phi_{\varepsilon}) \jm^{\varepsilon} (\phi_{\varepsilon}, t) .
  \label{eq:newmagtexteps}
\end{equation}
Computing as in {\tmname{Step 3}}, we have that
\begin{align}
  \mathcal{E}_{{\jom}}( \ju_{\varepsilon}) & \eqs 
  \int_0^{2 \pi} \int_{t \in I} \frac{| \jm^{\varepsilon}_{\bot}
  (\phi_{\varepsilon}, t) |^2}{\mathfrak{h}^2_1 (t)} + \frac{|
  \partial_t \jm^{\varepsilon} (\phi_{\varepsilon}, t)
  |^2}{\mathfrak{h}^2_2 (t)} \mathd \phi \sqrt{\mathfrak{g} (t)} \mathd
  t  \label{eq:equality0eps}\\
  &   \qquad  \qquad + \int_0^{2 \pi} \int_{t \in I} g (
  \jm^{\varepsilon} (\phi_{\varepsilon}, t) \cdot \tmmathbf{a}
  (\phi_{\varepsilon}, t) ) \mathd \phi \sqrt{\mathfrak{g} (t)} \mathd t
  \\
  & \eqs \int_0^{2 \pi}
  \Phi_{\mathcal{E}}^{\varepsilon} (\phi_{\varepsilon}) \mathd \phi \\
  & \leqslant  \int_0^{2 \pi} \Phi_{\mathcal{E}}^{\varepsilon} (\phi) \mathd
  \phi \\
  & \eqs  \int_{t \in I} \int_0^{2 \pi} \frac{|
  \jm^{\varepsilon}_{\bot} (\phi, t) |^2}{\mathfrak{h}^2_1 (t)} +
  \frac{| \partial_t \jm^{\varepsilon} (\phi, t)
  |^2}{\mathfrak{h}^2_2 (t)} \sqrt{\mathfrak{g} (t)} \mathd \phi \mathd
  t \nonumber\\
  &   \qquad \qquad  + \int_0^{2 \pi} \int_{t \in I} g (
  \jm^{\varepsilon} (\phi, t) \cdot \tmmathbf{a} (\phi, t) ) \mathd \phi
  \sqrt{\mathfrak{g} (t)} \mathd t  \label{eq:equality1eps}\\
  & \leqslant  \int_{t \in I} \int_0^{2
  \pi}  \frac{| \partial_{\phi} \jm^{\varepsilon}_{\bot} (\phi, t)
  |^2}{\mathfrak{h}^2_1 (t)} + \frac{| \partial_t \jm^{\varepsilon}
  (\phi, t) |^2}{\mathfrak{h}^2_2 (t)} \mathd \phi \sqrt{\mathfrak{g}
  (t)} \mathd t \\
  &  \qquad  \qquad + \int_{t \in I} \JW^2 (t) | \langle
  \jm^{\varepsilon}_{\bot} \rangle (t) |^2 \sqrt{\mathfrak{g} (t)}
  \mathd t \nonumber\\
  &   \qquad \qquad \qquad \qquad  + \int_0^{2 \pi} \int_{t \in I} g (
  \jm^{\varepsilon} (\phi, t) \cdot \tmmathbf{a} (\phi, t) ) \mathd \phi
  \sqrt{\mathfrak{g} (t)} \mathd t  \label{eq:equality2eps}\\
  & \leqslant \mathcal{E}_{{\jom}}( \jm^{\varepsilon}
  ) .  \label{eq:equality3eps}
\end{align}
Now, we observe that since $\mathcal{E}_{{\jom}} (
\jm^{\varepsilon} )$ is uniformly bounded, so is
$\mathcal{E}_{{\jom}} ( \ju_{\varepsilon} )$.
Therefore, there exists $\ju \in H^1 ( \Ssurf, \Tsurf )$ such that
$\ju_{\varepsilon} \rightharpoonup \ju$ weakly in $H^1 ( \Ssurf, \Tsurf
)$ and, by compactness, up to a subsequence, one also has that
$\ju_{\varepsilon} \to \ju$ strongly in $L^2 ( \Ssurf, \Tsurf
)$ and $\ju_{\varepsilon} \to \ju$ a.e. in $\Ssurf$. But also,
the family $( \Rm (\phi) \ju_{\varepsilon} (\phi, t) = \Rm
(\phi_{\varepsilon}) \jm^{\varepsilon} (\phi_{\varepsilon}, t)
)_{\varepsilon}$ is bounded in $H^1 ( I, \Tsurf )$ and,
therefore, there exists $\ju_{\ast} \in H^1 ( I, \Tsurf )$ such
that $\Rm (\phi_{\varepsilon}) \jm^{\varepsilon} (\phi_{\varepsilon}, t)
\rightharpoonup \ju_{\ast}$ weakly in $H^1 ( I, \Tsurf )$. By
compactness, up to a subsequence, one also has that $\Rm (\phi_{\varepsilon})
\jm^{\varepsilon} (\phi_{\varepsilon}, t) \to \ju_{\ast}$ \ strongly
in $L^2 ( I, \Tsurf )$ and $\Rm (\phi_{\varepsilon})
\jm^{\varepsilon} (\phi_{\varepsilon}, t) \to \ju_{\ast}$ a.e. in $I$
and, due to the independence of the family $( \Rm (\phi_{\varepsilon})
\jm^{\varepsilon} (\phi_{\varepsilon}, t) )_{\varepsilon}$ on the
$\phi$-variable, the same convergence relations also hold, respectively,
weakly in $H^1 ( \Ssurf, \Tsurf )$, strongly in $L^2( \Ssurf,
\Tsurf )$ and a.e. in $\Ssurf$.

A first consequence of the previous convergence relations is that the limit
vector field $\ju$ is axially symmetric. Indeed, we know that
\begin{align}
  &\Rm (\phi) \ju_{\varepsilon} (\phi, t) \to  \Rm (\phi) \ju
  (\phi, t)   \text{ a.e. in } \Ssurf, \\
  &\Rm (\phi) \ju_{\varepsilon} (\phi, t)  \to  \ju_{\ast} (t) 
   \text{ a.e. in } \Ssurf, 
\end{align}
and, therefore, also that $\ju (\phi, t) = \Rm^{\top} (\phi) \ju_{\ast} (t)$.
This shows that $\ju$ is axially symmetric.

Also, by the weak lower semicontinuity of the norm, the strong convergence in
$L^2( \Ssurf, \Tsurf )$, from \eqref{eq:equality0eps} and
\eqref{eq:equality3eps} we conclude that
\begin{equation}
  \mathcal{E}_{{\jom}}( \ju ) \leqslant
  \liminf_{\varepsilon \to 0} \mathcal{E}_{{\jom}}
  ( \ju_{\varepsilon} ) \leqslant \lim_{\varepsilon \to 0}
  \text{\eqref{eq:equality2eps}} \; \leqslant \lim_{\varepsilon \to 0}
  \mathcal{E}_{{\jom}} ( \jm^{\varepsilon} )
  =\mathcal{E}_{{\jom}} ( \jm ) .
\end{equation}
By the minimality of $\jm$, all previous inequalities are actually equalities,
and we can pass to the limit under the integral sign in
\eqref{eq:equality2eps}. Hence, we are back to the hypotheses necessary to resume the argument from {\tmname{Step 4}} and conclude the proof. Indeed, {\tmname{Step 4}} starts from the assumptions that if $\jm$ is a minimizer of
$\mathcal{E}_{{\jom}}$ then there exists an axially symmetric
vector field $\ju \in H^1( \Ssurf, \Tsurf )$, built from $\jm$,
such that $\mathcal{E}_{\jom}( \ju ) =\mathcal{E}_{\jom} (
\jm )$, and this step does not rely on any smoothness assumption on $\jm$.

\section{Further results and applications}\label{sec:furthresappls}
\noindent In this last section, we want to emphasize a particular case of
Theorem~\ref{thm:main1}. Formally, when $\jom (\xi) \assign \jlam \in \RR_+$
is constant and $\jlam \to + \infty$, the energy functional
$\mathcal{E}_{\jlam}$ converges, in the sense of $\Gamma$-convergence, to the
energy functional
\begin{equation}
  \mathcal{E} ( \jm ) \assign \mathcal{D}_{\Ssurf} ( \jm
  ) +\mathcal{A}_{\Ssurf} ( \jm ) \label{eq:eneunconstrained}
\end{equation}
subject to the constraint that $\langle \jm \rangle_{\Ssurf_{\xi}}
\times \ee_3 = 0$ for every $\xi \in \Ssurf$. It is convenient for us to
introduce the following notion.

\begin{definition} \label{def:axiallynullav}
  We say that a vector field $\jm \in H^1 ( \Ssurf, \Tsurf )$ is
  {\tmem{axially null-average}}, along the $\ee_3$-axis, if $\langle \jm
  \rangle_{\Ssurf_{\xi}} \times \ee_3 = 0$ for every $\xi \in \Ssurf$.
  Equivalently, given that $|\langle \jm
  \rangle_{\Ssurf_{\xi}} \times \ee_3 |^2 = | \langle
  \jm_{\bot} \rangle_{\Ssurf_{\xi}} |^2$, the vector field $\jm$
  is axially null-average when $\langle \jm_{\bot}
  \rangle_{\Ssurf_{\xi}} = 0$ for every $\xi \in \Ssurf$. In local
  coordinates, $\jm$ is axially null-average if, and only if, for every $t\in I$ there holds that
  \begin{equation}
    \langle \jm (\cdot, t) \rangle \assign \frac{1}{2 \pi}
    \int_0^{2 \pi} \jm_{\bot} (\theta, t) \mathd \theta \eqs 0 \, ,
     \label{eq:condprojnull}
  \end{equation}
  where $\jm_{\bot} (\theta, t) \assign ( \jm_{\bot} \circ \xip ) (\theta,
  t)$.
\end{definition}

\begin{remark}
  \label{eq:rmkaxsymm}Vector fields with line symmetry are axially
  null-average (i.e., satisfies \eqref{eq:condprojnull}). In particular, so
  are axially symmetric and antisymmetric configurations. Indeed, if $\jm$ is
  axially symmetric with respect to the $\ee_3$-axis then, in local
  coordinates, we have that
  \begin{equation}
    \jm (\phi, t) = \Rm^{\top} (\phi) \tmmathbf{\alpha} (t) \quad \forall
    (\phi, t) \in \RR \times I
  \end{equation}
  for some profile $\tmmathbf{\alpha} \in H^1 ( I, \Tsurf )$.
  Hence, $\langle \jm (\cdot, t) \rangle = ( \tmmathbf{\alpha}
  (t) \cdot \ee_3 ) \ee_3$ for every $t \in I$, and this implies that
  $\langle \jm_{\bot} \rangle_{\Ssurf_{\xi}} = 0$ for every $\xi
  \in \Ssurf$. A similar argument shows that, more generally, vector fields
  with line symmetry are axially null-average.
  
  Note that the class of axially null-average vector fields is not directly
  related to the class of null-average configurations in $H^1 ( \Ssurf,
  \Tsurf )$. Even if $\jm$ is $t$-invariant, \eqref{eq:condprojnull}
  does not imply that $\jm$ is null-average, but only that its projection
  $\jm_{\bot}$ is null-average.
\end{remark}

The main idea behind the proofs of Theorems \ref{thm:main0} and
\ref{thm:main1} is geometric and relatively intuitive. Yet it brings up many
interesting consequences that are not easy to report in a systematic way. We
want to mention the last one. As a byproduct of our analysis, we get that if
$\jm$ is a minimizer of $\mathcal{E}$ in $H^1 ( \Ssurf, \Tsurf )$,
which, a posteriori, turns out to be axially null-average, then any minimizer
of $\mathcal{E}$ is of the form \eqref{eq:genexprmperp} and, if in addition,
$\Tsurf$ is never-flat, then any minimizer of $\mathcal{E}$ has line symmetry.
Indeed, one can repeat verbatim the argument used to prove
Theorems~\ref{thm:main0} and \ref{thm:main1} to obtain the following result.

\begin{theorem}
  \label{thm:main3}For given surfaces of revolution $\Ssurf, \Tsurf \subseteq
  \RR^3$, we consider the energy functional defined for every $\jm \in H^1
  ( \Ssurf, \Tsurf )$ by
  {\opt}cf.{\tmem{~\eqref{eq:eneunconstrained}}}{\cpt}
  \[ \mathcal{E}  ( \jm ) \assign \int_{\Ssurf} | \grad \jm
     (\xi) |^2 \mathd \xi + \int_{\Ssurf} g ( \jm (\xi) \cdot
     \tmmathbf{a} (\xi) ) \mathd \xi \]
  with $g : \RR \to \RR_+$ a Lipschitz function and $\tmmathbf{a}:
  \Ssurf \to \RR^3$ a Lipschitz vector field. The following assertions
  hold:
  \begin{enumerateroman}
    \item If the vector field $\tmmathbf{a}$ is axially symmetric or
    antisymmetric, and $\jm \in H^1 ( \Ssurf, \Tsurf )$ is an
    axially null-average minimizer of $\mathcal{E}$, then $\jm$ is of the form
    {\opt}cf.~{\tmem{\eqref{eq:genexprmperp}}}{\cpt}
    \[ \jm (\phi, t) =\tmmathbf{\alpha}_{\bot} (t) \cos \phi
       +\tmmathbf{\beta}_{\bot} (t) \sin \phi + \eta (t) \ee_3 \]
    for some $\tmmathbf{\alpha}_{\bot}, \tmmathbf{\beta}_{\bot} \in H^1 (
    I, \RR^2 )$, $\eta \in H^1 ( I, \RR)$, subject to the
    constraint of the resulting $\jm$ being $\Tsurf$-valued.
    
    \item If $\Tsurf$ is never flat and $\tmmathbf{a}$ is axially symmetric,
    and $\jm \in H^1 ( \Ssurf, \Tsurf )$ is an axially null-average
    minimizer of $\mathcal{E}$, then either every minimizer is axially
    symmetric, or there also coexist axially antisymmetric minimizers. If
    $\Tsurf$ is never flat and $\tmmathbf{a}$ is axially antisymmetric, then
    either every minimizer is axially antisymmetric, or there also coexist
    axially symmetric minimizers. In any case, if the vector field
    $\tmmathbf{a}$ is axially symmetric or antisymmetric, then any axially
    null-average minimizer $\jm \in H^1 ( \Ssurf, \Tsurf )$ of
    $\mathcal{E}$ has line symmetry and is of the form
    \[ \jm (\phi, t) =\tmmathbf{\alpha}_{\bot} (t) \cos \phi
       +\tmmathbf{\beta}_{\bot} (t) \sin \phi + \eta (t) \ee_3 . \]
    for some $\tmmathbf{\alpha}_{\bot}, \tmmathbf{\beta}_{\bot} \in H^1 (
    I, \RR^2 )$, $\eta \in H^1 ( I, \RR )$, subject to the
    orthogonality conditions $\tmmathbf{\beta}_{\bot} (t) \cdot
    \tmmathbf{\alpha}_{\bot} (t) \equiv 0$ and $| \tmmathbf{\beta}_{\bot} (t)
    | \equiv | \tmmathbf{\alpha}_{\bot} (t) |$, and to the conditions of the
    resulting $\jm$ being $\Tsurf$-valued.
  \end{enumerateroman}
\end{theorem}

\begin{remark}
We want to highlight the main message contained in Theorem~\ref{thm:main3}: while it is always the case that axially symmetric configurations are axially null average, minimality allows for the converse implication.
\end{remark}
\begin{remark}
  Note that in the statement of Theorem~\ref{thm:main3}
  there is no more reference to the weight $\jom$ and to the circular integral
  weight $\JW$.  
\end{remark}

\begin{remark}
  \label{rmk:important} We are presenting Theorem~\ref{thm:main3} in a concise form that particularizes only the main assertions contained in Theorems~\ref{thm:main0}~and~\ref{thm:main1} to the current context of axially null-average configurations. However, \emph{every} claim in Theorems~\ref{thm:main0}~and~\ref{thm:main1} transposes to the current setting. In particular, the results about the existence of axially symmetric (and axially antisymmetric) minimizers specified in Theorem~\ref{thm:main0}.{\tmem{i}}~and~\ref{thm:main0}.{\tmem{ii}}  still apply to axially null-average minimizers.
\end{remark}

To better explain the interest in Theorem~\ref{thm:main3}, one can imagine the
following scenario, which will be explained through a concrete example soon,
in which by writing down the Euler--Lagrange equations associated with
$\mathcal{E}$ in \eqref{eq:eneunconstrained}, one infers that every stationary
point of $\mathcal{E}$ has to be {\tmem{axially null-average}}. In other
words, if by any means one can prove that for some specific choices of
$\Ssurf$ and $\Tsurf$ minimizers of $\mathcal{E}$ are necessarily axially
null-average, then the statements of Theorems \ref{thm:main0} and
\ref{thm:main1} still hold {\tmem{without}} any assumption on the circular
integral weight $\JW$.

\begin{example}
  Let us consider the simplest case in which $\Ssurf \assign \mathbb{A}$ is
  the annulus of $\RR^2$ of inner radius $t_1$ and outer radius $t_2$,
  $\Tsurf$ is the whole space $\RR^2$, and $g ( \jm (\xi) \cdot
  \tmmathbf{a} (\xi) ) = \kappa ( \jm (\xi) \cdot \ee_3 )^2$
  for some $\kappa \in \RR$. The functional $\mathcal{E} $ has to be minimized
  among all possible vector fields in $H^1 ( \mathbb{A}, \RR^2 )$
  that satisfy the boundary condition $\jm =\tmmathbf{b}_1$ on $t_1 \Stwo^1$
  and $\jm =\tmmathbf{b}_2$ on $t_2 \Stwo^1$ for axially symmetric vector
  fields $\tmmathbf{b}_1, \tmmathbf{b}_2$. We want to understand the symmetry
  properties of the minimizers. For that, we observe that any minimizer of
  $\mathcal{E}$ has to satisfy the Euler--Lagrange equations
  \begin{equation}
    - \mathLaplace \jm + \kappa ( \jm \cdot \ee_3 ) \ee_3  \eqs 0
    \quad \text{in } \mathbb{A}, \label{eq:pdeexample}
  \end{equation}
  under the prescribed boundary conditions $\jm =\tmmathbf{b}_1$ on $t_1
  \Stwo^1$ and $\jm =\tmmathbf{b}_2$ on $t_2 \Stwo^1$. By the previous
  relation, it follows that any minimizer of $\mathcal{E}$ is such that $-
  \mathLaplace \jm_{\bot} \eqs 0$ in $\mathbb{A}$. Expressing $\mathbb{A}$ in
  local coordinates through the classical polar map $\xip (\phi, t) =
  \Rm^{\top} (\phi) \jgam (t)$, $\jgam (t) \assign t \ee_1$, case in which the
  metric coefficients are $\mathfrak{h}^2_1 (t) = t^2$, $\mathfrak{h}^2_2 (t)
  = 1$, and then integrating in the $\phi$-variable, we conclude that any
  minimizer of $\mathcal{E}$ is such that
  \begin{equation}
    - \partial_t ( t \partial_t \langle \jm_{\bot}\rangle (t)
    ) = 0 \quad \text{in } [t_1, t_2]\, ,
  \end{equation}
  and subject to the boundary conditions $\langle \jm_{\bot}
  \rangle (t_1) = \langle \jm_{\bot} \rangle (t_2) = 0$. But
  the only solution of this boundary value problem is the zero solution, i.e.,
  the solution $\langle \jm_{\bot} \rangle (t) = 0$ for every $t
  \in [t_1, t_2]$. It follows that any solution of the boundary value problem
  \eqref{eq:pdeexample} is {\tmem{axially null-average}} and, therefore, the
  conclusions of Theorem~\ref{thm:main3} apply. But as pointed out in
  Remark~\ref{rmk:important}, also the conclusions of Theorem~\ref{thm:main0}
  apply from which we infer that there always exist axially symmetric
  solutions of the boundary value problem \eqref{eq:pdeexample}.
\end{example}

\section*{Acknowledgements}
\noindent{\sc G.DiF.} acknowledges support from the Austrian Science Fund (FWF)through the project {\emph{Analysis and Modeling of Magnetic
    Skyrmions}} (grant P34609).    {\sc G.DiF.} also
thanks TU Wien and MedUni Wien for their support and hospitality.
{\sc V.V.S.} acknowledges support by Leverhulme grant RPG-2018-438. 
{\sc A.Z.} has been partially supported by the Basque Government through the BERC 2022-2025 program and by the Spanish State Research Agency through BCAM Severo Ochoa excellence accreditation SEV-2017-0718 and through project PID2020-114189RB-I00 funded by Agencia Estatal de Investigacion (PID2020-114189RB-I00 / AEI / 10.13039/501100011033). A.Z. was also partially supported  by a grant of the Ministry of Research, Innovation and Digitization, CNCS - UEFISCDI, project number PN-III-P4-PCE-2021-0921, within PNCDI III.

{\sc G.DiF.}  and {\sc V.V.S.} would like to thank the Max
Planck Institute for Mathematics in the Sciences in Leipzig for
support and hospitality. {\sc G.DiF.} and {\sc A.Z.} thank the Hausdorff Research Institute for Mathematics in Bonn for its hospitality during the Trimester Program \emph{Mathematics for Complex Materials} funded by the Deutsche
Forschungsgemeinschaft (DFG, German Research Foundation) under Germany
Excellence Strategy -- EXC-2047/1 -- 390685813.

\bibliographystyle{acm}
\bibliography{literature}

\end{document}